\documentclass[12pt]{article}
\usepackage{amsthm,amssymb,amsmath}

\newcommand{\pP}{\mathbb P}
\newcommand{\Q}{\mathbb Q}
\newcommand{\R}{\mathbb R}
\newcommand{\Z}{\mathbb Z}

\newcommand{\sB}{\mathcal B}
\newcommand{\sD}{\mathcal D}
\newcommand{\sP}{\mathcal P}

\newcommand{\oB}{\overline{B}}

\newcommand{\oV}{\overline{V}}

\newcommand{\oX}{\overline{X}}

\newcommand{\ep}{\varepsilon}

\newcommand{\rddown}[1]{\left\lfloor{#1}\right\rfloor} 

\DeclareMathOperator{\cent}{center}
\DeclareMathOperator{\LCS}{LCS}

\DeclareMathOperator{\cNE}{\overline{NE}}
\DeclareMathOperator{\Supp}{Supp}

\theoremstyle{plain}
 \newtheorem{add}{Addendum}
 \newtheorem{cor}{Corollary}
 \newtheorem{lemma}{Lemma}
 \newtheorem*{mlemma}{Main Lemma}
 \newtheorem{prop}{Proposition}
 \newtheorem*{rth}{Revised Reduction}
 \newtheorem*{rind}{Revised Induction}
 \newtheorem{theorem}{Theorem}

\theoremstyle{definition}
 \newtheorem{df}{Definition}

\theoremstyle{remark}
 \newtheorem*{exa}{Example}
 \newtheorem*{que}{Question}
 \newtheorem{cau}{Caution}

\newcommand{\step}[1]{\paragraph{{\sc Step #1:}}}

\title{LETTERS OF A BI-RATIONALIST\\
{VII.~Ordered termination}}

\date{June~5th/July~24th, 2006 --
Серебряная свадьба; Moscow}

\author{V.V. Shokurov
\thanks{Partially supported
by NSF grant DMS-0400832.}
}

\begin{document}

\maketitle

\begin{abstract}
To construct a resulting model in LMMP, it is sufficient to prove
existence of log flips and their termination for certain
sequences. We prove that LMMP in dimension $d-1$ and termination
of terminal log flips in dimension $d$ imply, for any log pair of
dimension $d$, existence of a {\em resulting} log model: a
strictly log minimal model or a strictly log terminal Mori log
fibration, and imply existence of log flips in dimension $d+1$. As
consequence, we prove existence of a resulting model of $4$-fold
log pairs, existence of log flips in dimension $5$, and
Geography of log models in dimension $4$.
\end{abstract}

\ \ \ \ \ \ \ \ Числа не знаем, бо кончать не маем ...

\ \ \ \ \ \ \ \ Из письма зпорожских казаков турецькому султану.

\bigskip
The main purpose of this note to show that under certain inductive
assumptions and termination of terminal flips
we can construct for any log pair either
its log minimal model or a Mori log fibration.
This amounts to a weaker form of
the Log Minimal Model Program (LMMP) in which
termination of any sequence log flips is
replace by termination of some sequences.
This idea polishes the reduction to pl flips
\cite[4.5 and Section~6]{Sh92} \cite{Sh00}, and
appears recently also in \cite{AHK} (cf. Definition~\ref{Hter} below).
It looks that, for most of applications of LMMP,
it is sufficient.
Up to dimension $4$, we can omit the inductive assumptions and
the termination.
However the results may not be exaggerated.
They only demonstrate that progress in LMMP
hinge upon that of in termination.

LMMP means that of \cite[Section~5]{sh96b}.
By a terminal log flip we mean a log flip or
a divisorial contraction having only terminal
points in the flipping or exceptional locus respectively, that is,
the minimal log discrepancy (mld) for those points is $>1$.
Terminal termination: any sequence of
terminal log flips is finite.
Usually we apply this termination
to extremal and $\Q$-factorial contractions at least
near the flipping locus (cf. Caution~\ref{qfac} in
the proof of Theorem~\ref{fin}).

We work over a base field $k$ of characteristic $0$;
in some instances, $k$ is algebraically closed
or we need slightly change meaning of extremal curve
(see Definition~\ref{ecurve}).
We use standard facts and notation of LMMP, as
in \cite{ISh} \cite{KMM} \cite{Sh00}.
In particular, we use standard abbreviations:
{\em lt\/} for log terminal;
{\em dlt\/} for divisorially lt;
{\em klt\/} for Kawamata lt;
{\em lc\/} for log canonical;
{\em wlc\/} for weakly log canonical.
Let us briefly recall also some terminology
and notation:
birational {\em rational\/} $1$-contraction
does not blow up any divisor \cite[p.~84]{Sh00};
$\sD_B$ is the $\R$-vector space of Weil
$\R$-divisors $D$ supported in $\Supp B$
with the maximal absolute value norm
$\|D\|$ \cite[p.~134]{Sh00};
an FT (Fano type) relative variety means
a relative variety $X/Z$, for which there exists
an $\R$-boundary $B$,
such that $(X/Z,B)$ is a klt Fano pair;
for the cone of curves, contractions on FT $X/Z$
and for more details, see
\cite[Section~2]{PSh};
$K=K_X$ denotes a canonical divisor on a variety $X$;
$\LCS(X,B)$ is the subvariety of nonklt points.

\begin{theorem} \label{mod}
We assume LMMP in dimension $d-1$ and
termination of terminal log flips in dimension $d$.
Then any pair $(X/Z,B)$ of dimension $d$ with an $\R$-boundary $B$
has a {\em resulting} log model.
More precisely, $(X/Z,B)$ has either a strictly log minimal model, or
a strictly log terminal model with a Mori log fibration.

A pair $(X/Z,B)$ has a log minimal model
if and only if its numerical log Kodaira dimension
is nonnegative \cite[p.~263]{sh96b}, or equivalently
for lc $(X,B)$, $K+B$ is pseudo-effective.
\end{theorem}

Actually, LMMP also can be replaced by
terminal termination in dimension $\le d-1$
(cf. Corollary~\ref{exfl} and its proof).

\begin{que}
However,
if $B$ is a $\Q$-boundary, whether it is sufficient
LMMP with $\Q$-boundaries and
terminal termination with $\Q$-boundaries?
\end{que}

\begin{add} \label{termod}
If a starting pair $(X/Z,B)$ of dimension $d$
is strictly log terminal
or dlt a resulting model can be constructed by
a terminated sequence of (extremal) log flips.
\end{add}

Perhaps, for a more general starting pair, e.g.,
lc $(X/Z,B)$,
this part of LMMP also works: resulting models exist
and Geography of log models holds
in dimension $d$.
Similarly, under the same assumptions,
directed klt flops terminate  in dimension $d$.
In what follows, we  mention
some of those results, without the assumptions
for $4$-folds; we can add to this
termination of $4$-fold klt directed flops.

\begin{cor} \label{closed}
Under the assumptions of Theorem~\ref{mod},
the nonnegative numerical log Kodaira dimension is
a closed condition with respect to boundaries,
or equivalently, the birational existence of
a Mori log fibration is open.
\end{cor}

To obtain that result for usual log Kodaira
dimension, one needs semiampleness \cite[Conjecture~2.6]{sh96b}.

\begin{proof}
By definition the numerical dimension is
defined for a wlc model \cite[Proposition~2.4]{sh96b}, and
the wlc property is closed.
Moreover, for a limiting boundary of
such varieties a resulting model is wlc.
Otherwise by Theorem~\ref{mod} and Addendum~\ref{termod} we get
a Mori log fibrations on a birational model
that is open with respect to a boundary.
\end{proof}

\begin{cor}[Cf. {\cite[Theorem~3.4]{AHK}}]
Any $4$-fold pair $(X/Z,B)$ with an $\R$-boundary $B$
has a resulting model.
The closed and open properties of Corollary~\ref{closed}
hold for $4$-fold pairs.
\end{cor}

Note that semiampleness \cite[Conjecture~2.6]{sh96b},
special termination \cite[Example~8]{Sh04} and
\cite[Theorem~2.15]{AHK} improve the last result:
any sequence of log flips for a $4$-fold lc pair $(X/Z,B)$
with pseudo-effective $K+B$ terminates (see a remark
before Definition~\ref{Hter}).
Actually, semiampleness can be replaced by the effective
property of $(X/Z,B)$ in the sense of Birkar.
In particular, this gives termination of log flips if
$(X/Z,B)$ has the numerical log Kodaira dimension $0$
(not only big \cite[Corollary~3.6]{AHK}).

\begin{proof}
Immediate by LMMP up to dimension $3$ \cite[Theorem~5.2]{sh96b}
and terminal termination up to dimension $4$
\cite[Example~9 and Lemma~2]{Sh04}.
\end{proof}

\begin{rth}
LMMP in dimension $d-1$,
termination of terminal log flips in
dimension $d$ and existence of
pl flips in dimension $d+1$
imply existence of log flips
in dimension $d+1$.
\end{rth}

We recall that the flipping contraction is
assumed to be extremal and small,
the flipping variety is assumed to be $\Q$-factorial.

\begin{rind}
LMMP in dimension $d-1$ and
termination of terminal log flips in
dimension $d$
imply existence of pl flips
in dimension $d+1$.

Actually, it is sufficient to
assume the termination for {\em birational pairs\/} $(X/Z,B)$,
that is, $X\to Z$ is a birational contraction.
\end{rind}

The way we prove it below \cite{Sh00} \cite{HMc},
we can prove that any restricted
divisorial algebra on the reduced component
of a pl contraction in dimension $d+1$ is
finitely generated.
Similarly we can establish existence
of pl flips in dimension $n$ with
core dimension $n-s\le d$ \cite[1.1]{Sh00}.
So, we get a more terminal version of the main theorem
in \cite{HMc}.

\begin{cor} \label{exfl}
LMMP in dimension $d-1$ and
termination of terminal log flips in
dimension $d$
imply existence of log flips
in dimension $d+1$.
\end{cor}

As in Revised Reduction the flipping contraction is
assumed to be extremal and small,
the flipping variety is assumed to be $\Q$-factorial.
For $d$-dimensional log flips, those conditions
can be omitted.

\begin{proof}
Immediate by Revised Reduction and Induction.

Notice that in dimension $d$ we can obtain more
general log flips as in usual Reduction \cite[Theorem~1.2]{Sh00}.
\end{proof}

\begin{cor}
Log flips as in Corollary~\ref{exfl}
exist up to dimension $5$.
\end{cor}

A different proof see in \cite[Theorem~4.3]{AHK}.
Note that for $4$-fold log flips we need
only terminal termination in dimension $3$
\cite[Theorem~3.5]{ISh}.

\begin{proof}
Immediate by LMMP up to dimension $3$ \cite[Theorem~5.2]{sh96b}
and terminal termination up to dimension $4$
\cite[Example~9 and Lemma~2]{Sh04}.
\end{proof}

Other applications are as follows.

\begin{cor}[see {\cite[Theorem~6.20]{sh96b}
\cite[Conjecture~2.10]{ISh}}] \label{geog}
Geography conjecture holds for klt relative
$4$-folds.
This gives a bounded termination for
$D$-flips on any relative FT $4$-fold.
\end{cor}

This can be done for any lc relative pairs
after a slight generalization of Proposition~\ref{dcc}
and its addenda, corollaries to lc singularities
instead of dlt (see a remark after
Proposition~\ref{dcc}).

\begin{proof}
Geography can be obtain from Theorems~\ref{mod} and \ref{fin} below
(cf. the proof of Corollary~\ref{exfl}).

In the case of FT varieties this gives
a universal bound for $D$-termination.

A detailed treatment will be done elsewhere.

\end{proof}

\begin{cor}[see {\cite[Theorem~3.33]{Sh00}}]
Existence of Zariski decomposition for
relative {\em FT\/} $4$-folds.
In particular, any divisorial algebra for
any $\Q$-divisor is finitely generated on
such a variety.
\end{cor}

\begin{proof}
Immediate by Corollary~\ref{geog}.
\end{proof}

\begin{cor}
If a complete algebraic space of dimension $4$ has
only klt singularities and no rational
curves over an algebraic closure of
the base field, then it is projective.
\end{cor}

\begin{proof}
We can use methods of \cite{Sh95} and
Addendum~\ref{termod}.
\end{proof}

\begin{df} \label{ecurve}
An irreducible curve $C$ on $X/Z$ is called {\em extremal\/}
if it generates an extremal ray $R=\R_+[C]$ of
the Kleiman-Mori cone $\cNE(X/Z)$, and has
the minimal degree for this ray
(with respect to any ample divisor).
We also suppose that $R$ is contractible.

If the base field is not algebraically closed
then the contraction can be not extremal over its algebraic
closure, and
we take, for an extremal curve, that of
a (partial) extremal subcontraction, or
a sum of conjugations of such a curve divided by the number
of curves in the orbit.
\end{df}

If $(X/Z,B)$ is a dlt log pair with a boundary $B$,
such that $K+B$ has index $m$ then,
for any extremal curve $C/Z$,
$$
(K+B,C)\in \{\frac{n}{m}\mid n\in\Z, \text{ and }
n\ge -2dm\},
$$
where $d=\dim X$.
Immediate by the anticanonical boundedness
\cite[Theorem]{sh96a}.
In particular, $(K+C,B)\ge 1/m$ if $(K+B,C)>0$
(cf. \cite[Lemma~6.19]{sh96b}).

We can generalize these results for $\R$-boundaries.

\begin{prop} \label{dcc}
Let $(X/Z,B)$ be a lc pair with an $\R$-boundary $B$.
Then there exists
a finite set of real positive numbers $r_i$, and
a positive integer $m$ such that for any
extremal curve $C/Z$, for which $(X,B)$ is dlt
near the generic point of $C$,
$$
(K+B,C)\in \{\sum \frac{r_i n_i}{m}\mid
n_i \in\Z,\text{ and } n_i\ge -2dm\},
$$
where $d=\dim X$.
\end{prop}

If $(X,B)$ is dlt everywhere we can take any
extremal curve.
It is expected that actually we can relax dlt to lc
in the proposition, in its addenda and the corollaries:
LMMP is sufficient for this in dimension $d$
\cite[Conjecture and Heuristic Arguments]{sh96b}.
More precisely, it is sufficient existence of a strictly log minimal
model over any lc pair.
For this it is enough log flips and special termination
in dimension $d$ that follows from
LMMP in dimension $d-1$ by Corollary~\ref{exfl},
\cite{HMc} and \cite[Theorem~2.3]{Sh00}.
In addition, the lc property is better than the dlt one:
the former is closed (see Example and the proof of
Corollary~\ref{stabray} below).
However in our applications, dlt is the only assumption
what we need.

\begin{add}\label{adcc}
The numbers $r_i$ and $m,d$ depend on a pair $(X/Z,B)$
but they are the same after a (generalized) log flop
outside $\LCS(X,B)$, that is, only in curves $C$
with $C\cap\LCS(X,B)=\emptyset$.
\end{add}

To determine these numbers we use the following fact.

\begin{lemma}\label{dec}
Under the assumptions of Proposition~\ref{dcc},
there exists a decomposition
$B=\sum r_i B_i$ where $r_i$ are positive real numbers,
and $B_i$ are (Weil) $\Q$-boundaries such that:
\begin{description}

\item{\rm (1)}
$\sum r_i=1$;

\item{\rm (2)}
each $\Supp B_i\subseteq \Supp B$;

\item{\rm (3)}
each $K+B_i$ is a $\Q$-Cartier divisor which
has the trivial intersection $(K+B_i,C)=0$
for any curve $C/Z$ with $(K+B,C)=0$;

\item{\rm (4)}
$$
K+B=\sum r_i(K+B_i);
$$

\item{\rm (5)}
each $(X,B_i)$ is lc, $\LCS(X,B_i)\subseteq\LCS(X,B)$, and
$(X,B_i)$ is dlt in the locus where so does $(X,B)$.
\end{description}
\end{lemma}

The last assumption is meaningful because
there exists the maximal dlt set in $X$ and it is open:
the complement to the closure of log canonical
centers which are not dlt.

\begin{proof}
The main problem here is possible real multiplicities
of $B$.
Property (4) immediate by (1).
To satisfy (2-3) we consider an affine
$\R$-space of $\R$-divisors
$$
\sD_B^0=\{D\mid\Supp D\subseteq \Supp B,
\text{ and } K+D \text{ satisfies the intersection
condition of (3)}\}.
$$
The last means that $K+D$ is $\R$-Cartier, and
$(K+D,C)=0$ for any curve $C/Z$ with $(K+B,C)=0$.
This space is actually finite dimensional
and defined over $\Q$.
More precisely, this is a finite dimensional $\R$-space$/\Q$ in
the finite dimensional $\R$-space $\sD_B$
of $\R$-Weil divisors supported
in $\Supp B$.
Note for this, that the $\R$-Cartier condition
gives a linear subspace over $\Q$, and
any canonical divisor $K$ is integral.
Thus the condition for $K+D$ to be $\R$-Cartier
gives a finite dimensional affine subspace over $\Q$.
Each condition $(K+D,C)=0$ is also rational
because each intersection $(K+B_i,C)=m_i$ is
rational.
Finally, any $D\in\sD_B^0$ is an affine (weighted) linear
combination of $\Q$-Cartier divisors $K+B_i$
with $B_i$ (so far not necessarily boundaries)
supported in $\Supp B$.
Note that $B\in \sD_B^0$.

On the other hand, the $\R$-boundaries $D\in \sD_B^0$,
with lc $(X,D)$,
form a convex closed rational polyhedron \cite[1.3.2]{Sh92}, and
$B$ belongs to this polyhedron.
Any small perturbation inside
of the polyhedron
preserves the klt property outside $\LCS(X,B)$  and
the dlt property of (5).
Therefore, $K+B$ has a required decomposition (cf. Step~1
in the proof of Corollary~\ref{stabray} below).
\end{proof}

\begin{proof}[Proof of Proposition~\ref{dcc}]
The numbers $r_i$ were introduced in Lemma~\ref{dec}.
The positive integer $m$ is
an index for all divisors $K+B_i$, that is,
each $m(K+B_i)$ is Cartier.
By (5) of Lemma~\ref{dec} and the anticanonical boundedness
\cite[Theorem]{sh96a} each $(K+B_i,C)\ge -2d$.
Thus
$$
(K+B,C)=\sum r_i(K+B_i,C)=\sum r_i\frac{n_i}{m},
$$
where $n_i \in\Z,\text{ and } n_i\ge -2dm$.

Finally, we prove the addendum.
For simplicity, we consider only
usual log flops (some remarks about
more general flops see below).
By those we mean birational rational $1$-contractions
$X\dashrightarrow X'/Z$ which and its inverse are
indetermined only in curves $C/Z$ with
$(K+B,C)=0$.
Note that the decomposition $B=\sum r_i B_i$
with all its properties is preserved under
log flops in those curves \cite[Definition~3.2]{ISh}.
The same holds for the dimension $d$ and
index $m$ \cite[2.9.1]{Sh83}.
The space $\sD_B^0$ and the boundary polyhedron are also preserved (for
small flops), or surjective on the corresponding
space, the polyhedron for any log flop.
For (5), it is enough that the log flop is
outside $\LCS(X,B)$.

A generalized log flop is a crepant modification
which can blow up some exceptional divisors with
log discrepancies $\le 1$ and $>0$
outside $\LCS(X,B)$, that is, with
centers not in $\LCS(X,B)$.
Kawamata says: such a modification
of log pairs is a {\em log $K+B$-equivalence\/}, and
the pairs are {\em log $K+B$-equivalent\/}.
For example, it can be a crepant blowup, or
its composition with subsequent log flops.
Under certain assumption (see Lemma~\ref{blow} below), such
a blowup exists according to
the finiteness of exceptional divisors with
log discrepancies $\le 1$ outside $\LCS(X,B)$.
Then the $\Q$-Cartier property of $K+B$ is enough on
a blowup, or even on a log resolution.
This condition is preserved under log flops even
generalized because the intersection numbers
can be computed on any common resolution
by the projection formula.

\end{proof}

\begin{cor}[on an interval; {cf. \cite[Lemma~6.19]{sh96b}}]
\label{int}
Let $(X/Z,B)$ be a lc pair with an $\R$-boundary $B$.
Then there exists a real number $\hbar>0$
such that for any extremal curve $C$, for which $(X,B)$ is dlt
near the generic point of $C$,
either

$(K+B,C)\ge \hbar$, or

$(K+B,C)\le 0$.
\end{cor}

That is, the intersection numbers $(K+B,C)$
do not belong to the interval $(0,\hbar)$.

\begin{add}
The number $\hbar$ depends on a pair $(X/Z,B)$
but it is the same after a (generalized) log flop
outside $\LCS(X,B)$.
\end{add}

Note that extremal curves may not
be preserved under flops!

\begin{proof}
By Proposition~\ref{dcc}, the intersection
numbers $(K+B,C)=\sum r_in_i/m$ with extremal curves $C/Z$
under the assumptions of the proposition,
satisfy the dcc.
Moreover, for any real number $A$ the set
of those numbers $\le A$ is finite.
Thus
$$
\hbar=\min \{\sum \frac{r_i n_i}{m}>0\mid
n_i \in\Z,\text{ and } n_i\ge -2dm\}
$$
is a required positive number.

By Addendum~\ref{adcc} we can take the same $\hbar$
after any log flop outside $\LCS(X,B)$.
\end{proof}

\begin{exa}
Let $L_i,i=1,2,3$, be $3$ distinct lines
in the plane $\pP^2$ passing through
a point $P$.
Then, for $F=L_1+L_2+L_3$,
$$
\sP=\{D\in \sD_F\mid (\pP^2,D)
\text{ is a lc pair with an $\R$-boundary \/} D\}
$$
is a convex closed rational polyhedron.
The face of boundaries $D=\sum b_i L_i$
with $\sum b_i=2$ gives nondlt pairs $(\pP^2,D)$
with $3$ exceptions: $(\pP^2,F-L_i), i=1,2,3$.
Thus the dlt property is not closed and not convex.
However, the dlt property holds
exactly in $(\pP^2\setminus P,D)$
for any interior point $D$ of the face.
\end{exa}

\begin{cor}[stability of extremal rays] \label{stabray}
Let $(X/Z,B)$ be a lc pair with
an $\R$-boundary $B$, and $F$ be a reduced divisor on $X$.
Then there exists a real number $\ep>0$
such that for any other $\R$-boundary $B'\in \sD_F$ and
any extremal contractible ray $R\subset\cNE(X/Z)$ such that:
\begin{description}

 \item{\rm (1)}
$\|B'-B\|<\ep$;

 \item{\rm (2)}
$K+B'$ is $\R$-Cartier, and $(K+B',R)<0$;

 \item{\rm (3)}
for some extremal curve $C$ in $R$, the pair $(X,B)$ is
dlt near the generic point of $C$, and so does $(X,B')$
possibly for a different extremal curve;

\end{description}
the seminegativity $(K+B,R)\le 0$ holds.

If $(X,B+E)$ is a dlt pair with
an effective $\R$-divisor $E$ and $\Supp (B+E)=F$
than we can omit condition (3).

\end{cor}

\begin{add} \label{adstab1}
For fixed $B'$, a log flop of $(X/Z,B)$
outside of $\LCS(X,B)$ in
any $R$ as in the corollary preserves
$\ep$ in direction $B'$, that is,
the stability holds again
for any $D$ in the segment $[B_Y,B_Y']$
on flopped $(Y/Z,B_Y)$, where $B_Y$ and $B_Y'$
denote the birational transforms on $Y$
of corresponding boundaries $B$ and $B'$ from $X$.
More precisely,  $D\in \sD_{F_Y}$, is an $\R$-boundary,
the log flop projects $\sD_F$ onto
$\sD_{F_Y}$ (some components of
$F$ are contracted),
where $F_Y$ denotes the birational transform of $F$ on $Y$;
\begin{description}

 \item{\rm (1)}
$\|B_Y-D\|<\ep$.

 \item{\rm (2)}
For any extremal contractible ray $R\subset\cNE(Y/Z)$
such that $(K_Y+D,R)<0$ and

 \item{\rm (3)}
for some extremal curve $C$ in $R$, the pair $(Y,B_Y)$ is
dlt near the generic point of $C$, and so does $(Y,B_Y')$
possibly for a different extremal curve;

\end{description}
$(K_Y+D,R)\le 0$ holds.
\end{add}

\begin{cau}
In other directions, log flops can spoil
singularities of $(Y,D)$.
\end{cau}

\begin{add} \label{adstab2}
The same applies to any log flop which is
a composition of such log flops as in Addendum~\ref{adstab1}, and
of log flops with one weaker condition:
\begin{description}

 \item{\rm (2)'}
$(K+B,R)=(K+B',R)=0$.

\end{description}
\end{add}

\begin{proof}
The main idea is that the dlt property is {\em conical\/}.

 \step1 {\em Choice of $\delta$.\/}
We can suppose that $B\in \sD_F$.
There exists a real number $\delta>0$ such that:
\begin{description}

 \item{\rm (3)'\/}
for any $\R$-boundary $B'\in \sD_F$ with $\R$-Cartier $K+B'$ and
any $D$ in the {\em open\/} ray $\overrightarrow{BB'}$
with $\|D-B\|\le \delta$,
$D$ is an $\R$-boundary,
the pair $(X,D)$ is dlt (exactly) in the locus where
so does $(X,B')$; in particular, $(X,D)$ satisfies
(3) for the same extremal curve as $(X,B')$.

\end{description}
By \cite[1.3.2]{Sh92}
$$
\sP_B=\{D\in \sD_F\mid (X,D)
\text{ is a lc pair with an $\R$-boundary \/} D\}
$$
is a convex closed (rational) polyhedral cone in
some $\delta$-neighborhood of $B$.
Unfortunately, a similar set for the dlt property,
instead of the lc one, can be not closed
(see Example above).
However the dlt property of $(X,D)$ holds
in the same maximal open subset of $X$
for all $D$ in the interior of every face of $\sP_B$.
Indeed, according to the linear behavior of discrepancies
with respect to $D$, all $D$ in the interior have
the same support on $X$ and the same
log canonical centers in $X$.
Thus a dlt resolution of $(X,D)$ over
the maximal dlt open subset in $X$ gives that of
for any other divisor in the interior.
Now we obtain (3)' from the dlt property in
the interior of a (minimal) face of $\sP_B$ with
$D=B'\in \sP_B$.

By monotonicity and stability \cite[1.3.3--4]{Sh92}, under
the last assumption in the corollary, the dlt property
is open and closed in $\sP_B$ near $B$.
Moreover, in the definition of $\sP_B$,
we can replace the lc property by
the $\R$-Cartier one near $B$.
Then we do not need (3).
Indeed, each prime component of $E$
does not passes through any log canonical center of $(X,B)$.

 \step2 {\em Required\/}
$$
\ep=\frac{\delta}{N+1},
$$
where $N$ is any positive number $\ge 2d/\hbar$,
and $d=\dim X$.

Indeed, if $C$ is an extremal curve in $R$
as in (2-3) of the corollary,
and $(K+B,C)>0$.
Then by Corollary~\ref{int},
$(K+B,C)\ge \hbar$, and
$$
(B-B',C)=(K+B,C)-(K+B',C)>\hbar.
$$
Hence
$$
(K+B'+N(B'-B),C)=(K+B',C)+N(B'-B,C)<-N\hbar\le -2d,
$$
that contradicts the anticanonical boundedness
\cite[Theorem]{sh96a}.
Indeed, by our choice of $\ep$ and (3)',
$D=B'+N(B'-B)=B+(N+1)(B'-B)$ is an $\R$-boundary,
$\|B-D\|=(N+1)\|B'-B\|<(N+1)\ep\le \delta$,
and $D$ satifies (3)'.
Thus $(K+B,R)\le 0$.

 \step4 {\em Addenda.\/}
In both addenda we consider a log flop in $R$,
that is, $(K+B,R)=0$.
Then we can replace $B'$ by a {\rm maximal\/}
$B'$ in direction $B'$, that is, by
a (possibly) new boundary $B'$ on
the ray $\overrightarrow{BB'}$ such that
$\|B'-B\|=\ep$, or infinitely close to $\ep$;
in the former case we need slightly decrease $\ep$.
(Since we use the maximal norm, actually in
most directions we can take $B'$ on
a larger Euclidian distance.)
Then any $B'$ and any $D'\in [B,B']$
will satisfy the properties (1-3) of
the corollary and (3)' of its proof.
In particular, (2) holds because
$(K+B,R)=0$.
The same applies to (2)' in Addenda~\ref{adstab2}.

The property (1) of $D$ in
the first addendum follow
almost by definition.
The distance in (1) is $<\ep$ and
less than the length of $[B_Y,B_Y']$
which can be shorter even $0$ when $B_Y'=B_Y$.
If the birational rational $1$-contraction $X\dashrightarrow Y$
contracts divisors, $F_Y$ is smaller than $F$, and
we need to change $\ep$ on $\ep_Y=\|B_Y'\|$,
which is the maximal norm with respect to
the (noncontracted) prime components of $F_Y$.
(Since $F$ has finitely many components,
$\ep_Y$ stabilizes after finitely many log flops.)
Property (3)' on $Y$ in direction $D_Y'$ with
$\delta_Y=(N+1)\ep_Y$
can be obtained from the fact that
log flips and log flops preserve or improve
log singularities.
Note that log flops outside
$\LCS(X,B)$ preserves klt, dlt singularities and
$\LCS(X,B)$ itself.
In Addendum~\ref{adstab1},
the log flop is a log flip with respect
to $K+D'$ if $D'\not =B$.
Such a log flip with respect to $K+D'$
improves singularities of $(X,D')$,
actually, (3)' holds after the log flop
for any $D$ in the segment $[B_Y,(N+1)B_Y']$;
$D'\to D_Y'=D$ is a surjective projection.
Thus by Addendum~\ref{adcc}, the constants
$\hbar,N,$ and $d$ are the same, $\ep_Y,\delta_Y$
as above,
and under the conditions (2-3) of the addendum
the required seminegativity holds.
This gives the same constants for log flops with $(K+D',R)=0$.

The second addendum follows by induction
on the composition.

\end{proof}

\begin{cor}[cf. {\cite[Corollary~6.18]{sh96b}}]
\label{sample}
Let $(X/Z,B)$ be a klt minimal model with
big $K+B/Z$.
Then it has a lc model, that is,
$K+B$ is semiample.
\end{cor}

\begin{proof}
This is well-known for $\Q$-divisor $B$ \cite[Theorem~2.1]{Sh83}
\cite[Theorem~3-1-1 and Remark~3-1-2]{KMM}.
By Lemma~\ref{dec} we have a decomposition
$B=\sum r_iB_i$ satisfying (1-5) of
the lemma.
Moreover, for sufficiently small real number $\ep>0$,
if $\|B_i-B\|<0$
we can suppose that $(X/Z,B_i)$ is klt  with big
$K+B_i/Z$, and by Corollary~\ref{stabray},
$K+B_i$ is nef$/Z$.
According to construction of the decomposition,
we can always find boundaries $B_i$ in
the $\ep$-neighborhood of $B$.
Therefore, each $K+B_i$ is semiample$/Z$, and
so does $K+B$.
Moreover, $(X/Z,B'=\sum r_i'B_i)$, with any
$0<r_i'\in\Q,\sum r_i'=1$,
is a klt minimal model with big $K+B'$
and with the same lc model as for $(X/Z,B)$;
all those models are {\em equivalent\/}
in the sense of \cite[Definition~6.1]{sh96b}.

\end{proof}

\begin{cor}[Stability of wlc models] \label{stabwlc}
Let $(X/Z,B)$ be a dlc wlc pair with an $\R$-boundary $B$,
$E$ be an effective $\R$-divisor, and $F$ be a reduced divisors
such that $(X,B+E)$ is a dlt pair and $\Supp (B+E)=F$.
Then there exists a real number $\ep>0$ such
that for any other $\R$-boundary $B'\in \sD_F$
with $\|B'-B\|<\ep$ the following statements
are equivalent:
\begin{description}

 \item{\rm (1)\/}
$H=B'-B$ is nef on any (irreducible) curve $C/Z$
with $(K+B,C)=0$;

 \item{\rm (2)\/}
for some real number $0<\delta<\ep/\|H\|$,
$(X/Z,B+\delta H)$ is a dlc wlc pair; and

 \item{\rm (3)\/}
for any real number $0<\delta<\ep/\|H\|$,
$(X/Z,B+\delta H)$ is a dlc wlc pair.

\end{description}
\end{cor}

\begin{add} \label{stabwlcad}
The models in (3) are equivalent.
\end{add}

\begin{proof}
We choose the same $\ep$ as in the proof
of Corollary~\ref{stabray}.
By our assumption the dlt property near $B$
coincides with the lc one: if, in Step~1 of the proof
of Corollary~\ref{stabray},
we replace the lc condition
by the dlc one, and even by $\R$-Cartier one,
we obtain the same cone
$\sP_B$ in the $\ep$-neighborhood of $B$.

(1)$\Rightarrow$(3):
The nef property in (1) includes the $\R$-Cartier
property of $H$.
By our choice of $\ep,\delta$ and $B'$,
$D=B+\delta H\in \sP_B$, and $D$ is an $\R$-boundary.
Thus by definition of $\sP_B$,
the pair $(X/Z,D)$ is a dlt pair.
If it is not wlc, then $K+D$ is not nef$/Z$, and
by \cite[Theorem~2]{Am}
there exists an extremal contractible ray $R\subset\cNE(X/Z)$
satisfying conditions (2-3) of Corollary~\ref{stabray} with $B'=D$.
Hence by the corollary $(K+B,R)\le 0$, and
so $(K+B,R)=0$ by the wlc property of $(X/Z,B)$.
But then $(H,R)<0$ that contradicts (1).

(3)$\Rightarrow$(2):
Immediate by assumptions.

(2)$\Rightarrow$(1):
Suppose that $(H,C)<0$ for some (irreducible) curve $C/Z$
with $(K+B,C)=0$.
Then $(K+B+\delta H,C)=\delta(H,C)<0$ that contradicts
the wlc property in (2).

The equivalence of the addendum follows
from the linear property of intersections.
If $(K+B+\delta H,C)=0$ for some $\delta$ in (3)
then $(K+B+\delta H,C)=0$ for any $\delta$ in (3).
Otherwise, $(K+B+\delta H,C)<0$ for some $\delta$ in (3).
Similarly, if $(K+B+\delta H,C)>0$ for some $\delta$ in (3)
then $(K+B+\delta H,C)>0$ for any $\delta$ in (3).

\end{proof}

\begin{lemma}[Convexity of equivalence] \label{conv}
If two wlc models are equivalent then all
resulting models between them exist and are wlc equivalent to
each of two models.
\end{lemma}

\begin{proof}
Suppose that wlc models $(X/Z,B)$ and $(X/Z,B')$
are equivalent; by definition we can
suppose that they have the same variety $X/Z$.
We verify that any model $(X/Z,B'')$ between
them, that is, for any $B''\in [B,B']$, is wlc and equivalent
to each of above models; again we take the same variety $X/Z$.
This gives also existence of a resulting model for $(X/Z,B'')$.

Indeed, for some real numbers $\alpha,\beta\ge0,
\alpha+\beta=1$, $B''=\alpha B+\beta B'$.
Let $C/Z$ be a curve with $(K+B,C)=0$.
Then $(K+B',C)=0$ because models $(X/Z,B)$ and
$(X/Z,B')$ are equivalent.
Hence by the linear property of intersection
$$
(K+B'',C)=\alpha(K+B,C)+\beta(K+B',C)=0.
$$
Similarly, if $(K+B,C)>0$, then
$(K+B',C)$ and $(K+B'',C)>0$.
Thus $(X/Z,B'')$ is a wlc model, equivalent
to $(X/Z,B)$ and to $(X/Z,B')$.
Notice also that $(X,B'')$ is lc by
(cf. \cite[1.3.2]{Sh92}) and
each log discrepancy $a(E,X,B'')=
\alpha a(E,X,B)+\beta a(E,X,B')$.

\end{proof}

The following concept is a formalization of
the well-known method from \cite[4.5]{Sh92} for
reduction of log flips to pl flips, and
it is extremely important in our proofs.
The same concept under a different name,
directed flips, appears in \cite{AHK}.
However, it is redundant there: in \cite[Theorem~3.4
and Corollaries~3.5-6]{AHK} any sequence of
log flips terminates.

\begin{df}[$H$-termination] \label{Hter}
Let $(X/Z,B)$ be a lc pair, and
$H$ be an $\R$-divisor.
A sequence of log flips (not necessarily extremal)
$$
(X_1=X/Z,B_1=B)\dashrightarrow
(X_2=X_1^+/Z,B_2=B_1^+)\dashrightarrow
\dots
$$
is called $H$-{\em ordered\/} if we can
associate a real number $\lambda_i>0$ with
each flip $X_i\dashrightarrow X_{i+1}/Z$
such that:
\begin{description}

 \item{\rm (1)}
the numbers decrease:
$\lambda_1\ge \lambda_2\ge \dots$;

 \item{\rm (2)}
each flip $X_i\dashrightarrow X_{i+1}/Z$ is
a log flop with respect to $K_{X_i}+B_i+\lambda_i H_i$,
where $H_i$ is the birational image of $H$ on $X_i$;
and

 \item{\rm (3)}
each pair $(X_i/Z,B_i+\lambda_i H_i)$ is
a wlc model.

\end{description}
We say that the flip $X_i\dashrightarrow X_{i+1}/Z$
has the {\em level\/} $\lambda_i$ with respect to $H$.
So $H$-{\em termination\/} of a given sequence of
$H$-ordered log flips means that it
terminates, that is, finite.
\end{df}

It is clear that termination
of any sequence of nontrivial log flips implies
$H$-termination for any its $H$-ordering.
On the other hand,
$H$-termination for a sequence of log flips is sufficient
for its termination and allows
to construct a resulting model.

\begin{prop} \label{resmod}
For any klt starting model $(X/Z,B+\lambda_1 H)$,
$H$-termination implies existence of
a resulting model for $(X/Z,B)$, in particular,
termination of the corresponding log flips of $(X/Z,B)$.

Moreover, for a given starting model,
an $H$-ordered sequence of log flips exists
if the log flips exist in dimension $d=\dim X$.
\end{prop}

\begin{proof}
Let $X_n\dashrightarrow X_{n+1}/Z$ be
the last flip of level $\lambda_n$.
By definition
$(X_{n+1}/Z,B_{n+1}+\lambda_n H_{n+1})$ is
a wlc model, and $\lambda=\lambda_n>0$.
If $K_{X_{n+1}}+B_{n+1}$ is nef$/Z$, then
$(X_{n+1}/Z,B_{n+1})$ is a wlc model
of $(X/Z,B)$, a resulting model.

Otherwise $K_{X_{n+1}}+B_{n+1}$ is not nef$/Z$.
By induction we can suppose that
$(X_{n+1}/Z,B_{n+1}+\lambda_n H_{n+1})$ is klt
(see below the proof of existence of $H$-ordered flips).
Thus by Corollary~\ref{stabray}
with $(X,B)=(X_{n+1},B_{n+1}+\lambda_n H_{n+1})$
either $\cNE(X_{n+1}/Z)$ has an extremal ray $R$
such that $(K_{X_{n+1}}+B_{n+1}+\lambda_n H_{n+1},R)=0$,
$(H_{n+1},R)>0$, and $(K_{X_{n+1}}+B_{n+1},R)<0$,
or by Corollary~\ref{stabwlc}
there exists $0< \lambda_{n+1}<\lambda_n$
such that $(X_{n+1}/Z,B_{n+1}+\lambda_{n+1} H_{n+1})$
is a wlc model.
In the former case, by our assumptions,
$R$ gives a Mori log fibration $X_{n+1}\to Y/Z$,
a resulting model with the boundary $B_{n+1}$
for $(X/Z,B)$, because $H_{n+1}$ is numerically ample$/Y$
and $\lambda_n>0$.
In the latter case, we proceed as follows.
We can assume that $\lambda_{n+1}$
is minimal in our construction, that is,
with nef $K_{X_{n+1}}+B_{n+1}+\lambda_{n+1} H_{n+1}/Z$;
the klt prperty is preserved by monotonicity \cite[1.3.3]{Sh92}.
Then we obtain a Mori log fibration as in the former case.

Now we explain how to extend the sequence
of log flips if $X_n\dashrightarrow X_{n+1}$ is
not a last one.
If the above contraction for $R$ is birational,
it has a log flip $X_{n+1}\dashrightarrow X_{n+2}=X_{n+1}^+/Z$
(possibly a divisorial contraction).
It is a log flop of the wlc model
$(X_{n+1}/Z,B_{n+1}+\lambda_n H_{n+1})$, and
thus satisfies (1-3) of Definition~\ref{Hter}
with $\lambda_{n+1}=\lambda_n$.
Note that it preserves the klt property of
the pair.
Otherwise we consider a similar construction
for minimal $\lambda_{n+1}<\lambda_n$ as above.
Again it extends $H$-ordered sequence, if
$R$ corresponds to a birational contraction.
Since $(X_{n+1},B_{n+1})$ is lc,
$(X_{n+1},B_{n+1}+\lambda_n H_{n+1})$
is klt, and $0<\lambda_{n+1}<\lambda_n$,
by monotonicity \cite[1.3.3]{Sh92}
the pair $(X_{n+1},B_{n+1}+\lambda_{n+1} H_{n+1})$
is klt.
Any flop preserves the klt property.
Hence $(X_{n+2},B_{n+2}+\lambda_{n+1} H_{n+2})$
is also klt that complets induction.

Usually we include into log flips
the divisorial contractions.
Thus we can consider not only small
modifications as log flops, or use the fact
that, after finitely many log flips, all
the next ones are small and so do the log flops.
\end{proof}

It is easy to give an example of
a sequence of log flips which cannot be
$H$-ordered at least for some divisor $H$.
Take two disjoint birational contractions,
one positive, another negative with respect to $H$.

In what follows,
all isomorphism of models, e.g., local ones, are induced by their
birational isomorphisms.

\begin{theorem} \label{fin}
We assume LMMP in dimension $d-1$ and
termination of terminal log flips in dimension $d$.
Let $(X_i/Z,B_i)$ be a sequence of $d$-dimensional
dlt wlc pairs which converge to a dlt pair $(X/Z,B)$ in
the following sense:
\begin{description}

 \item{\rm (1)}
each $X_i$ is isomorphic to $X$ (and between
themselves) in codimension $1$;
all divisors $B_i$ and $B$ are {\em finitely\/}
supported, that is, there exists a reduced
divisor $F$ such that $B$ and each $B_i\in \sD_F$;

 \item{\rm (2)}
each $X_i$ is isomorphic to $X$ near $\LCS(X,B)=S=\rddown{B}$
and $\LCS(X_i,B_i)=\LCS(X,B)$:
there exists neighborhoods $U_i$ of $\LCS(X_i,B_i)$
and $V_i$ of $\LCS(X,B)$ which are isomorphic and
identified under the birational isomorphism of (1);

 \item{\rm (3)}
there exist finitely many
prime b-divisors (exceptional and nonexceptional)
$D_j$ {\em outside\/} (possibly not
{\em disjoint\/} from) $\LCS(X,B)$, that is,
$\cent_X D_j\not\subseteq \LCS(X,B)$, and
which contain all positive {\em log codiscrepancies\/}
$b(D_j,X_i,B_i)=1-a(D_j,X_i,B_i)$ outside $\LCS(X,B)$,
that is, if $\cent_{X_i}D_j\not\subseteq\LCS(X,B)$
and $b(D_j,X_i,B_i)> 0$ for some $i$,
then $D_j$ is one these b-divisors;

 \item{\rm (4)}
there exists a limit of {\rm b-$\R$-divisors\/}
$\oB_i=S+\sum b(D_j,X_i,B_i)D_j$:
$$
\oB_{\lim}=S+\sum b_jD_j=S+\sum\lim_{i\to\infty}
b(D_j,X_i,B_i)D_j;
$$
and

 \item{\rm (5)}
$B=B_{\lim}=S+\sum b_j D_j$, where
the summation runs only nonexceptional $D_j$ on $X$, and
$\oB\ge\oB_{\lim}$, where
$\oB=S+\sum b(D_j,X,B)D_j$ is the crepant b-subboundary for $(X,B)$
extended in the b-divisors $D_j$.

\end{description}
Then the sequence is finite in the {\em  model sense\/},
that is, the set of equivalence classes of
models $(X_i/Z,B_i)$ is finite.
\end{theorem}

Notice that actually (3) implies existence
of a finite support in (1).

\begin{cor} \label{equivmod}
We assume LMMP in dimension $d-1$ and
termination of terminal log flips in dimension $d$.
Let $(X_i/Z,B_i)$ be a sequence of $d$-dimentional
dlt log pairs such that:
\begin{description}

 \item{\rm (1)\/}
each model $(X_i/Z,B_i)$ is wlc;

 \item{\rm (2)\/}
the models are isomorphic in codimension $1$;
and isomorphic near $\LCS(X_i,B_i)$;

 \item{\rm (3)\/}
for some $\R$-boundaries $B$ and $B'$,
each $B_i\in [B,B']$, and the models are ordered
in the segment: $B_i=B+\lambda_i H, H=B'-B,
\lambda_1\ge \lambda_2 \ge \dots, \lambda_i\in (0,1]$;
and

 \item{\rm (4)}
for some $i$, $(X_i/Z,B)$ and $(X_i/Z,B')$
are dlt log pairs
with $\LCS(X_i,B)=\LCS(X_i,B')=\LCS(X_i,B_i)$.

\end{description}
Then the models {\em stabilizes\/}:
the models are equivalent for $i\gg 0$.
\end{cor}

Note that (3) is meaningful because the Weil divisors
on each model are the same by
the first statement of (2).

\begin{proof}
\

 \step1 {\em Nonequivalence of models.\/}
By Lemma~\ref{conv} we can suppose that
numbers $\lambda_i$ form an infinite sequence,
$\lambda_0=\lim_{i\to\infty} \lambda_i$, and
the models $(X_i,B_i)$ are pairwise
nonequivalent.
Otherwise the stabilization holds.

 \step2 {\em Conditions of Theorem~\ref{fin} are
satisfied for an appropriate subsequence.\/}
We can suppose $i=1$ in our assumption (4).
Take $(X=X_1/Z,B:=B_{\lim})$ where
$B:=B_{\lim}=\lim_{i\to\infty}B_i=B+\lambda_0 H$, or
$\lambda_0=0$ for new $B$.
Conditions (1--2) of Theorem~\ref{fin} hold
by assumptions (2-4).
We can make condition (3) in Theorem~\ref{fin} taking
a subset of $D_j$ with $b(D_j,X,B)\ge 0$, and
$\cent_X D_j\not\subseteq \LCS(X,B)$, or
equivalently, $\cent_X D_j\not\in X\setminus \LCS(X,B)$;
if $D_j$ is also nonexceptional, it
is assumed that $D_j$ is supported in $\Supp B$ or in $\Supp B'$.
The set of $D_j$ is finite by \cite[Corollary~1.7]{sh96b}.
Then the condition holds for $(X,B_i)$ with
all $B_i$ sufficiently close to $B$ by assumptions (2) and (4),
by stability of the klt property and continuity
of log discrepancies with respect to the multiplicities in $D_j$,
where $B_i$ on $X$ is its birational transform from $X_i$.
Hence assumption (1) and monotonicity \cite[Lemma~2.4]{ISh}
imply (3) in Theorem~\ref{fin} for $i\gg 0$.

Up to convergency in (4) of Theorem~\ref{fin},
the conditions in (5) follow from construction and
monotonicity \cite[Lemma~2.4]{ISh}.
Indeed, $B=B_{\lim}$ by construction,
$b(D_j,X,B_i)\ge b(D_j,X_i,B_i)$
by the monotonicity and wlc property of $(X_i/Z,B_i)$.
Thus
$$
b(D_j,X,B)=b(D_j,X,B_{\lim})=
\lim_{i\to\infty} b(D_j,X,B_i)\ge
\lim_{i\to\infty} b(D_j,X_i,B_i)=b_j,
$$
and $\oB\ge\oB_{\lim}$
which gives (5) of Theorem~\ref{fin}.

If, for fixed $D_j$,
$b(D_j,X_i,B_i)$ is not bounded from below,
we can drop such $D_j$, and take
a subsequence with $\lim_{i\to\infty} b(D_j,X_i,B_i)=-\infty$.
In the bounded case, we have finally a convergent subsequence
in (4).

Now the stabilization follows from the finiteness
in Theorem~\ref{fin} that contradicts Step~1.
\end{proof}

\begin{lemma}[Canonical blowup] \label{blow}
We assume LMMP in dimension $d-1$ and
termination of terminal log flips in dimension $d$.
Let $(X,B)$ be a klt log pair of dimension $d$, and
$Z$ a closed subvariety of codimension $\ge 2$.
Then there exists a (unique under the algorithm
in the proof) crepant blowup $Y\to X$
such that
\begin{description}

 \item{\rm (1)}
$Y$ is isomorphic to $X$ over $X\setminus Z$;

 \item{\rm (2)}
$(Y,B_Y)$ is cn in codimension $\ge 2$ over $Z$;

 \item{\rm (3)}
if in addition, $K_Y+B$ with the birational
transform of boundary $B$ on $Y$ is ample,
a blow up is unique.

\end{description}
\end{lemma}

\begin{proof}[Proof-Construction]
\

 \step1
Consider a log resolution $(Y/X,B^+)$ with
a boundary $B^+=\sum b_j^+ D_j$ combined with
log codiscrepancies: $b_j^+=\max\{b(D_j,X,B),0\}$.
We can suppose that the prime  components of
$\Supp B^+$ are disjoint, and there exists
finitely many exceptional divisors $E/X$
with $b(E,X,B)\ge 0$, or equivalently
with $a(E,X,B)\le 1$, by the klt property
\cite[Corollary~1.7]{sh96b}.
Moreover, all $b_j^+<1$, and
$(Y,B^+)$ is terminal in codimension $\ge 2$.
We use slightly different boundary $B_Z\le B^+$ on $X$:
$0$ in all exceptional divisors$/Z$, and
$B^+$ elsewhere.

\step2
We apply LMMP to $(Y/X,B_Z)$.
The log flips exist by \cite[Theorem~1.1]{HMc} or induction
of Corollary~\ref{exfl}.
Each flip is terminal because we never
contract nonzero components $E$ of $B_Z$.
Indeed, this holds over $Z$ by
our assumptions because corresponding
boundary multiplicities are zero.
Otherwise we get a component $E$ with
$P=\cent_X E\not\subset Z$, and
$(X/X,B)$ near $P$ is the lc model
of $(Y/X,B_Z)$, even after the divisorial contraction.
Since this contraction decreases the codiscrepancy,
and increases the discrepancy in $E$ that
contradicts \cite[Lemma~2.4]{ISh}.
Thus termination holds by our assumptions.
Since the resulting model $(Y/X,B_Z)$ is birational$/X$,
it is a terminal in codimension $\ge 2$ and
a strictly log minimal model.

\step3
Using semiampleness in the big klt case we
obtain the lc model $(Y/X,B_Z)$;
the previous model $Y$ in Step~2 and this
model $Y$ are FT$/X$ (by Step~4 below).
The model satisfies (1) and (3).
Indeed, (1) follows
from the uniqueness of a lc model.
Then by construction $B_Z$ on $Y$ is
the birational transform of $B$, and $K_Y+B$
is ample and such a model is also unique.

\step4
By the negativity \cite[1.1]{Sh92}, for the crepant model
$(Y,B_Y)$ of $(X,B)$, the subboundary $B_Y$ is a boundary.
Still possible, $B_Y$ can have
noncanonical singularities of codimension $\ge 2$.
Then we can apply Steps 1--2 to $(Y,B_Y)$
with $Z_Y$ as union of all noncanonical centers$/Z$
of $(Y,B_Y)$.
This process is terminated and we
obtain finally a crepant
model $(Y/X,B_Y)$.
Indeed, each construction blows up at
least one exceptional divisor $E$ with
$a(E,X,B)> 0$ over a noncanonical center, and
there exists only
finitely many such divisors.
However, in general, $K_Y+B$ may be not ample
over $X$.
The above algorithm gives a unique blow up by (3).
\end{proof}

\begin{lemma}[$D$-flip] \label{flip}
We assume LMMP in dimension $d-1$ and
termination of terminal log flips in dimension $d$.
Let $(X,B)$ be a klt pair of dimension $d$, and
$D$ be a prime divisor on $X$ such that
$(X,B)$ is terminal in codimension $\ge 2$ at $D$, that is,
if $E$ is an exceptional divisor with
$a(E,X,B)\le 1$ then $\cent_X E\not\subseteq D$.
Then a $D$-flip of $X/X$ exists.
\end{lemma}

\begin{proof}[Proof-Construction]
The construction is quite standard (cf.
terminalization \cite[Theorem~6.5]{ISh} and
$\Q$-factorialization \cite[Lemma~7.8]{ISh}).
By the uniqueness of $D$-flips, they can be
constructed locally$/X$ \cite[Corollary~3.6]{Sh00}.

 \step1
As in Steps~1--2 of the proof of Lemma~\ref{blow}
with $Z=\emptyset$, that is, the starting $B_Z=B^+$,
we obtain a crepant blow up $(Y/X,B_Y)$ which
is terminal in codimension $\ge 2$ and strictly log minimal$/X$.
By our assumptions there are no exceptional
divisors$/D$.

 \step2
Let $D$ be its birational transform on $Y$.
Now we apply $D$-MMP to construct nef $D/X$.
For a sufficiently small real number $\ep>0$,
$(Y/X,B_Y+\ep D)$ is terminal in codimension $\ge 2$, and $D$-MMP is
LMMP for the pair.
Again the log flips exist.
Termination is terminal and holds:
$D$-flips do not contract any divisor.
Thus we can suppose that $D$ is nef$/X$.

 \step3
Contraction given by $D$ or $K_Y+B_Y+\ep D$
is the required model for $D$, $D$-flip.
The contraction exists as in Step~3
of the proof of Lemma~\ref{blow} because
$Y$ is FT$/X$.
Note that the model is small over $D$,
because $Y$ do not have exceptional divisors$/D$.
On the other fibers$/X$, $D$ is trivial,and
$Y$ is isomorphic to $X$ over $X\setminus D$.

\end{proof}

\begin{mlemma}
Let $(X/Z,B)$ be a dlt pair, $(X'/Z,B_{X'})$ be
its wlc model, isomorphic to $(X,B)$
near $\LCS(X,B)=\LCS(X',B_{X'})$, and
$X\to Y/Z$ be an (extremal) contraction negative with respect to
$K+B$.
Then the contraction is birational,
with the exceptional locus
disjoint from $\LCS(X,B)$, and
contracts only b-divisors $D$ with
$b(D,X,B)>b(D,X',B_{X'})$.
\end{mlemma}

A contracted b-divisor $D$ has
$\cent_X D$ in the contracted locus.

\begin{proof}
By our assumptions $B_{X'}^{\log}=B_{X'}$, that is,
$X\dashrightarrow X'$ is a rational $1$-contraction.
Thus by \cite[Proposition~2.5, (ii)]{ISh}
(cf. \cite[Proposition~2.4.1]{sh96b})
the contraction is not fibred.
In addition, for any irreducible curve $C/Z$
intersecting $\LCS(X,B)$,
$(K+B,C)\ge (K_{X'}+B_{X'},C')\ge 0$
where $C'$ is the birational image of $C$ on $X'$;
the latter is well defined by our assumptions.
Indeed, the log discrepancies for prime b-divisors
with centers near $\LCS(X,B)$, in particular,
for centers intersecting $\LCS(X,B)$, are
the same for $(X',B_{X'})$, and
by monotonicity \cite[Lemma~2.4]{ISh}
$b(D_i,X,B)\ge b(D_i,X',B_{X'})$
for the other b-divisors $D_i$.
This and the projection formula
for a common resolution of $X$ and $X'$
implies the inequality (cf. the proof of
\cite[Proposition~2.5, (i)]{ISh}).
In particular, the exceptional locus of $X/Y$
is disjoint from $\LCS(X,B)$.

Now for simplicity suppose that there exists
a log flip $(X^+/Y/Z,B^+)$ of $X/Y$: in our applications
we always have it.
(Otherwise one can use the last statement of
\cite[Lemma~2.4]{ISh}.)
Then $(X'/Z,B_{X'})$ is also
a wlc model of the divisorial contraction
or of the log flip
--- the basic fact of LMMP.
Thus $b(D,X,B)>b(D,X^+,B^+)\ge b(D,X',B_{X'})$, or
equivalently,
$a(D,X,B)<a(D,X^+,B^+)\le a(D,X',B_{X'})$
by monotonicities \cite[Lemmas~3.4 and 2.4]{ISh}.
\end{proof}

\begin{proof}[Proof of Theorem~\ref{fin}]
Taking a subsequence, we can suppose that
the models $(X_i/Z,B_i)$ pairwise are not
equivalent.
Then we need to verify that the sequence
is finite.

We care only about models outside $\LCS(X,B)$.
Near the $\LCS(X,B)$, the models are wlc by (2), and
we will keep this: assuming $F$ minimal in (1),
\begin{description}

 \item{\rm (6)}
in the proof below,
all models $(Y/Z,D)$ with an $\R$-boundary $D$ are isomorphic to
$X$ near $\LCS(Y,D)=\LCS(X,B)$,
$D\in \sD_F$ near $\LCS(X,B)$, and thus
there exists a real number $\ep>0$ such that
$(Y,D)$ is dlt near $\LCS(X,B)$ if
in addition $K_Y+D$ is $\R$-Cartier and $\|D-B'\|<\ep$.

\end{description}
Notice also that in our construction below
$B,B_i$ and similar $\R$-boundaries $D$ will have
the same reduced part: $\LCS(Y,D)=\rddown{D}=\LCS(X,B)=S$.

\step1 {\em Terminal limit\/}.
We construct a dlt model $(\oX/Z,B_{\oX})$
such that:
\begin{description}

 \item{\rm (7)}
$\oX,X$ and each $X_i$ are isomorphic
near $\LCS(X,B)=\LCS(\oX,B_{\oX})$;

 \item{\rm (8)}
$\oX\dashrightarrow X$ and each $\oX\dashrightarrow X_i$
are birational rational $1$-contractions;
$\oX$ blows up all $D_j$ of (3) with $b_j\ge 0$ and
with $\cent_X D_j\cap \LCS(X,B)=\emptyset$;

 \item{\rm (9)}
every $\R$-divisor $D$ on $\oX$ which
is $\R$-Cartier near $\LCS(X,B)$ is $\R$-Cartier
everywhere on $\oX$; in particular, each
divisor $D$ on $\oX$ with
$\Supp D\cap \LCS(X,B)=\emptyset$
is $\Q$-Cartier;

 \item{\rm (10)}
$B_{\oX}\ge \oB_{\lim}$ as b-divisors but divisors on $\oX$
(that is, for prime divisors on $\oX$); in particular,
for $D_j$ with nonnegative multiplicities $b_j$
and with $\cent_X D_j\cap \LCS(X,B)=\emptyset$ (see (8));
and

 \item{\rm (11)}
$(\oX,B_{\oX})$ is dlt, and terminal {\em completely\/}
outside $\LCS(X,B)$
in the following sense: $a(E,\oX,B_{\oX})>1$, or
equivalently $b(E,\oX,B_{\oX})<0$,
for each exceptional divisor with
$\cent_{\oX}E\cap \LCS(X,B)=\emptyset$.

\end{description}

By Lemma~\ref{blow} we can construct a slightly
weaker version with properties (7),
(8) for $b_j>0$ by (5), because
$b(D_j,X,B)\ge b_j>0$, and (10).
We apply the lemma to $(X,B)$ with the closed subvariety
which is the union of $\cent_{X} E$
for exceptional divisors $E$ with
$a(E,X,B)<1$, or
equivalently $b(E,X,B)>0$,
and with $\cent_X E\cap \LCS(X,B)=\emptyset$,
in particular,
of all $\cent_X D_j$ in (8) with $b_j>0$.

Since $\oX$ can be not $\Q$-factorial, we
slightly modify $\oX$ to be sufficiently
$\Q$-factorial and terminal.
To blow up the canonical centers completely outside
or {\em disjoint\/} from $\LCS(X,B)$,
we can use an increased divisor (boundary outside of $\LCS(X,B)$)
$B_{\oX}+\ep H$ where $H$ is a general
ample Cartier divisor passing through such
centers.
For a sufficiently small real number $\ep>0$,
the noncanonical centers of $(\oX,B_{\oX}+\ep H)$
are the only canonical $\cent_{\oX}E$'s
with $a(E,\oX,B_{\oX})=0$ and
$\cent_{\oX} E\cap \LCS(X,B)=\emptyset$ (cf. \cite[1.3.4]{Sh92}).
This gives (8) for $b_j=0$ by (5), (11) and preserves (7), (10).
To satisfy (9), it is enough to perform this for one divisor $D$
which is sufficiently general near $\LCS(X,B)$.
Indeed, by rationality of klt singularities,
the Weil $\R$-divisors modulo $\sim_\R/X\setminus \LCS(X,B)$ have
finitely many generators.
Since the $\R$-Cartier property defines
an $\R$-linear subspace over $\Q$ among $\R$-divisors,
we can suppose that such generators are
Cartier near $\LCS(X,B)$ and integral.
Adding ample divisors we can suppose that
they are prime and free near $\LCS(X,B)$, and
thus by (11) not passing through the canonical
(i.e., nonterminal) centers outside $\LCS(X,B)$
(even everywhere).
We can make $\Q$-Cartier each $D$ one by one.
According to Lemma~\ref{flip}, there exists
a small modification over $\oX$ such that
$D$ is $\Q$-Cartier on the modification ($D$-flip).
This gives (9) and concludes the step.
The dlt property of (11) near $\LCS(\oX,B_{\oX})$ by (6-7).

 \step2 {\em Limit of boundaries.\/}
For each $i$, let $B_i^+$ be an $\R$-boundary on $\oX$
with multiplicities $\max\{b(D,X_i,B_i),0\}$
in the prime divisors $D$ on $\oX$.
We can replace (10) by a more precise version:
\begin{description}

 \item{\rm (10)'}
$B_{\oX}= \oB_{\lim}^+=S+\sum b_j^+D_j,
b_j^+=\max\{b_j,0\},$ as b-boundaries, including
exceptional $D_j$ on $X$
with the nonnegative multiplicities $b_j=b_j^+$
and with $\cent_X D_j\cap \LCS(X,B)=\emptyset$
(cf. (3), (8) and (10) above),
and $0=b_j^+$ for all other $D_j$ with
$\cent_X D_j\cap \LCS(X,B)=\emptyset$
(and with $b_j<0$; and such
$D_j$ are possible); or equivalently,
$$
B_{\oX}=\lim_{i\to\infty} B_i^+.
$$

\end{description}
By monotonicity \cite[1.3.3]{Sh92} property (11) is preserved;
(6) and the other properties of $(\oX/Z,B_{\oX})$ so do.
By (9) the $\R$-Cartier property holds for all adjoint
divisors $K_{\oX}+B_i^+$ and $K_{\oX}+B_{\oX}$;
$B_{\oX}=B=B_{\lim}$ near $\LCS(X,B)$.

In addition, by (6-7), (9) and (11):
\begin{description}

 \item{\rm (11)'}
each $(\oX,B_i^+)$ is a dlt pair, terminal in the sense of (11);
$\LCS(\oX,B_i^+)=\LCS(\oX,B_{\oX})=\LCS(X,B)$; and
$B_i^+=B_i$ near $\LCS(X,B)$.

\end{description}
This is true by stability of terminal and klt singularities
(cf. \cite[1.3.4]{Sh92}) after taking a subsequence of models
$(X_i/Z,B_i)$ for all $i\gg 0$.
The last statement in (11)' allows to use for $B_i^+$
the properties of $B_i$ near $\LCS(X,B)$, e.g., (2).

 \step3 {\em Wlc terminal limit.\/}
We can suppose that $(\oX/Z,B_{\oX})$ is
a wlc model, terminal in the sense of (11).
Otherwise by (11)
there exists an extremal contraction
$\oX\to Y/Z$ negative with respect to $K_{\oX}+B_{\oX}$
\cite[Theorem~2]{Am}.
We claim that the contraction is birational,
does not contract components $D$ of $B_{\oX}$
with positive multiplicities, and does not
touch $\LCS(X,B)$, that is, an isomorphism
in a neighborhood of $\LCS(X,B)$.
Indeed, such a contraction is stable for
a small perturbation of the divisor $B_{\oX}$:
for any $\R$-boundary $B'\in\sD_{F+\sum D_j}$ sufficiently close
to $B_{\oX}$ and with $\R$-Cartier $K_{\oX}+B'$,
the contraction will be negative with respect
to $K_{\oX}+B'$.
By (10-11)', for all $i\gg 0$,
the contraction is negative with respect
to $K_{\oX}+B_i^+$.
By construction and definition $(X_i/Z,B_i)$ is
a wlc model of $(\oX/Z,B_i^+)$.
Notice that $B_i=(B_i^+)_{X_i}^{\log}$ because
$\oX\dashrightarrow X_i$ is a birational rational $1$-contraction by (1) and (8).
Therefore the contraction is not fibred by (6-7) and
Main Lemma, or
equivalently, it is birational.
The contraction is disjoint from $\LCS(X,B)$
by (6-7) and the same lemma.
Finally, the contraction does not contract
a prime divisor $D$ with positive multiplicity
in $B_{\oX}$ because by (10)' they are positive for
$B_i^+$ with $i\gg 0$, and $D=D_j$ with $b_j>0$.
It is impossible by Main Lemma again:
$b(D_j,\oX,B_i^+)=b(D_j,X_i,B_i)$.
Moreover, the contraction does not
contract $D_j$ with $b_j=0$.
Indeed, after such a contraction
$a+1=a(D_j,Y,B_Y)>a(D_j,\oX,B_{\oX})=1,a>0$, and
since boundaries $B_{Yi}^+$ for all $i\gg 0$ are
small perturbations of $B_Y$,
then $a(D_j,Y,B_{Yi}^+)\ge 1+a/2$ for all $i\gg 0$
where $B_Y,B_{Yi}^+$ are the images of
$B_{\oX},B_i^+$ respectively on $Y$.
Each $(X_i/Z,B_i)$ is also a wlc model
of $(Y/Z,B_{Yi}^+)$, and $a(D_j,X_i,B_i)\ge
a(D_j,Y,B_{Yi}^+)\ge 1+a/2$ by \cite[Lemma~2.4]{ISh},
or equivalently,
$b(D_j,X_i,B_i)\le -a/2$ and
$$
0=b_j=\lim_{i\to \infty}b(D_j,X_i,B_i)\le -a/2<0,
$$
a contradiction.

Therefore the contraction $\oX/Y$, if it is
divisorial, or otherwise its log flip, preserves
the properties (7-11) and (10-11)'
with the images of the corresponding
boundaries.
Note that by (9) the divisorial extremal contraction
blow down a $\Q$-Cartier divisor disjoint from
$\LCS(X,B)$, and, in particular, preserves (9).
The log flips exist by \cite[Theorem~1.1]{HMc} or induction
of Corollary~\ref{exfl}.
It preserves (9) by its extremal property
(cf. \cite[2.13.5]{Sh83}).
Of course, the property (11)' holds
after taking models for all $i\gg 0$.

Since each log flip is extremal by construction and
terminal by (11),
the flips terminate and we obtained
wlc $(\oX/Z,B_{\oX})$ with the required
terminal property.

\begin{cau} \label{qfac}
A log flip can be non-$\Q$-factorial, that is,
$\oX$ can be non-$\Q$-factorial.
However according to usual Reduction
\cite[Theorem~1.2]{Sh00} such a log flip exists
(cf. also the proof of Corollary~\ref{exfl} above).

Termination can be also non-$\Q$-factorial.
We can reduce it to usual $\Q$-factorial
terminal termination as in
the proof of special termination \cite[Theorem~4.8]{ISh}
taking a strictly log terminal blowup of $(\oX,B_{\oX})$;
for a dlt pair, a $\Q$-factorialization can be
constructed this way.
To construct such a model for any lc pair
in dimension $d$, it is sufficient
existence of $\Q$-factorial log flips and
special termination in this dimension.
Existence of terminal non-$\Q$-factorial log flips
follows from the same construction and
Corollary~\ref{sample} (cf. the proof of Lemma~\ref{flip}).
\end{cau}

 \step4 {\em Equivalence intervals.\/}
The intervals belon to the affine space $\sB$ of
$\R$-divisors on $\oX$ generated by
divisors $B_{\oX}$ and $B_i^+$.
It is a finite dimensional subspace in the linear space of
$\R$-divisors having the support in divisors $D_j$
and the birational transform of $F$ by (1), (10--11)' and (6--8).
In this affine space $B_{\oX}=\lim_{i\to\infty} B_i^+$.
Geography of log models \cite[Section~6]{sh96b} \cite[2.9]{ISh}
gives an expectation that near $B_{\oX}$, that is,
for boundaries in $\sB$ close to $B_{\oX}$,
there are only finitely many equivalent classes
of wlc models satisfying (6).
We prove it partially:
there exists a real number $\ep>0$ such
that, in each direction $B_i^+$, the wlc models
are equivalent in the interval of the length $\ep$.
Of course, we can assume that each $B_i^+\not=B_{\oX}$:
otherwise the model $(X_i/Z,B_i)$ is equivalent to
$(\oX/Z,B_{\oX})$.
Hence each direction is well-defined.
More precisely,
there exists an $\R$-boundary $B_i'\in\sB$ such that:
\begin{description}

 \item{\rm (12)}
$\|B_i'-B_{\oX}\|= \ep$;

 \item{\rm (13)}
$B_i^+\in (B_{\oX},B_i')$;
and

 \item{\rm (14)}
all wlc models in the interval $(B_{\oX},B_i')$
are equivalent:
for any $D\in(B_{\oX},B_i')$, $D\in \sB$ and
is an $\R$-boundary, $\rddown{D}=S$,
$(X_i/Z,D_{X_i})$ is a dlt wlc pair
equivalent to $(X_i/Z,B_i)$, and having
nonnegative codiscrepancies $b(D_j,X_i,D_{X_i})\ge 0$
only in the above $D_j$ where
$D_{X_i}$ is the image of $D$ on $X_i$;
$D$ as $B_i^+$ also satisfies (6);
$D_{X_i}$ as $B_i$ satisfies (2-3) and (6).

\end{description}
We can suppose that $\R$-boundaries $D\in\sB$ form actually
a cone with the vertex $B_{\oX}$ in
the $\ep$-neighborhood of $B_{\oX}$ (cf. \cite[1.3.2]{Sh92}).
To establish (14) we use $\ep$ from
Corollaries~\ref{stabray}, \ref{stabwlc} and
under the additional assumption:
\begin{description}

 \item{\rm (11)''}
for each $\R$-boundary $D\in\sB$ with $\|D-B_{\oX}\|\le\ep$,
with $\R$-Cartier $K_{\oX}+D$ and with $\rddown{D}=S$,
$(\oX,D)$ is a dlt pair satisfying (6) and
the terminal property of (11).

\end{description}
This follows from the stability of terminal and klt
singularities for small perturbation of $\R$-boundaries
\cite[1.3.4]{Sh92}.

Indeed, take $D=B_i'$ satisfying (12-13) and thus (11)''.
Then it is an $\R$-boundary and by (10)'
property (13) holds for all $i\gg 0$.
Now we apply LMMP to $(\oX/Z,B_i^+)$.
Again by (6) as in Step~3 above,
if $K_{\oX}+B_i^+$ is not nef
there exists an extremal negative
contraction with respect to $K_{\oX}+B_i^+$
\cite[Theorem~2]{Am}.
Its log flip exists by \cite[Theorem~1.1]{HMc} or induction
of Corollary~\ref{exfl}.
Termination holds as in Step 3.
(See also Caution~\ref{qfac} above.)
As in that step all extremal contractions
are birational, disjoint from $\LCS(X,B)$,
terminal, and do not contract
any prime component $D_j$ disjoint from
$\LCS(X,B)$ with $b(D_j,X_i,B_i)\ge 0$,
equal to multiplicities of $B_i^+$ in $D_j$;
other components of $B_i^+$ with multiplicity $0$
(even $D_j$) can be contracted.
According to terminal termination we
obtain a dlt wlc model $(\oX_i/Z,B_i^+)$, and
$(X_i/Z,B_i)$ is its wlc model.
They are equivalent; in particular,
$(\oX_i/Z,B_i^+)$ is a crepant model of $(X_i/Z,B_i)$
and $B_i^+=\oB_i$ on all divisors of $\oX_i$,
in particular, on $D_j$ blown up on $\oX_i$.
By Corollary~\ref{stabwlc} and its
Addendum~\ref{stabwlcad}
the same holds for any $D$ in the interval $(B_{\oX},B_i')$
with the corresponding model $(X_i/Z,D_{X_i})$.
For $E$ we take $eE$ where $E$ is the total support
for all $B_i^+-S$ and $0<e\ll 1$ is a real number.
The models are equivalent one each other
by Addendum~\ref{stabwlcad}.
Notice also that by Corollary~\ref{stabray}
the above log flips are
log flops with respect to $K_{\oX}+B_{\oX}$.
Thus in any direction,
$\ep$ is preserved by Addenda~\ref{adstab1}-\ref{adstab2}.
By (1) and (7--8) we can suppose that all $\Supp B_i^+$
are the same.
Then in the addenda
$\ep$ and $\delta$ are preserved for log flops
with respect to $K_{\oX}+B_{\oX}$:
the components of $\Supp B_i^+$ are not contracted
(see Step~4 in the proof of Corollary~\ref{stabray}).

 \step5 {\em $1$-dimensional case.\/}
If the real affine space $\sB$ has
dimension $1$, there are at most two
intervals $(B_{\oX},\pm B_i')$, and
at most two types of models.
For higher dimension we use:

 \step6 {\em Induction, or limit of equivalence intervals.\/}
The intervals $[B_{\oX},B_i']$ have
a convergent subsequence $\lim_{l\to\infty} B_{i_l}'=B'\in \sB$.
Otherwise we have finitely many intervals and
models as in Step~5.
The limit $B'$ is also an $\R$-boundary on $\oX$ and
by construction satisfies (11)''.
Now we cut the limit by an affine {\em rational\/} hyperplane:
there exists an affine rational hyperplane $\sB'\subset\sB$
such that it intersects $(B_{\oX},B')$ in $B_s$ and
the intervals $(B_{\oX},B_{i_l}')$ in $B_{sl}^+$.
The new boundaries $B_s,B_{sl}^+$ on $\oX$,
the images $B_{sl}$ of the latter ones on $X_{i_l}$
instead of $B_{\oX},B_{i_l}^+$, and $B_{i_l}$ respectively
satisfy the same properties (1--11), and (10--11)'
after taking a subsequence.
Thus the corresponding to $\sB$ space of divisors
is a subspace of $\sB'$.
(Actually we do not need $(X/Z,B)$ and corresponding properties;
$(\oX,B_{\oX}:=B_s)$ is sufficient.)
The properties (1--2), (6--7), (9), (11), and even (11)',
immediate by construction.
In (3) we can keep the same $D_j$ by (14).
If for some $D_j$, the set of codiscrepancies
$b(D_j,X_{i_l},B_{sl})$ is unbounded from
below, we take this subsequence and discard this $D_j$.
Therefore, we can find a subsequence satisfying (4).
Then (5), even (10), (10)', and (8) hold by construction.
For (10)', notice that $B_{sl}^+=\oB_{sl}$
for all $D_j$ with nonnegative multiplicities
in $\oB_{sl}$ by Step~4, and (14); and
extended by $0$ in other components of $\oB_{sl}$
exceptional on $X_{i_l}$ (cf. Step~2).

However $(\oX/Z,B_s)$ is not necessary wlc.
Therefore we apply again Step~3, etc.
This completes induction on
dimension of $\sB$.

\end{proof}

\begin{proof}[Proof-Construction of Theorem~\ref{mod}]
We construct strictly log terminal resulting
models $(X/Z,B_\lambda),B_\lambda=B^{\log}+\lambda H,$ for
some effective $\R$-divisor $H$ and $\lambda\in [0,1]$, and
find a real number $\lambda_0\in [0,1)$
such that $(X/Z,B_\lambda )$
are minimal for $\lambda \ge \lambda_0$, and
Mori log fibrations for $\lambda<\lambda_0$.
Thus we get a minimal model for $\lambda_0=0$, and
a Mori log fibration in all other cases.

 \step1
Using a Hironaka resolution we can suppose that
$(X/Z,B^{\log})$ is strictly log terminal.

 \step2
Then by special termination \cite[Theorem~2.3]{Sh00}
\cite[Corollary~4]{Sh04},
we can suppose that in any sequence of log flips
of $(X/Z,B^{\log})$ ($H$-ordered or not), the flips
are {\em nonspecial\/}, that is, do not
intersect $\LCS(X,B^{\log})$.
(For $\lambda_0=0$, this means that
$K+B^{\log}$ is nef on $\LCS(X,B^{\log})/Z$;
see Step~4 below.)

\step3
We can add a rather ample $\R$-boundary $H=\sum h_i H_i,h_i\not=0$,
with prime divisors $H_i$ such that
\begin{description}

 \item{\rm (1)}
$\rddown{H}=0$ and $\Supp H\cap \Supp B^{\log}=\emptyset$
in codimension $1$;

 \item{\rm (2)}
$(X/Z,B_1=B^{\log}+H)$ is
a strictly log minimal model;

 \item{\rm (3)}
prime components $H_i$ of $\Supp H$
generate the numerical classes of
all divisors$/Z$;
and

 \item{\rm (4)\/}
the multiplicities $h_i$ of $H$ are {\em independent\/}
over $\Q(B)$:
$$
\sum a_i h_i=a,\text{ all } a_i, a\in \Q(B)\Longrightarrow
\text{ all } a_i=0, a=0,
$$
where $\Q(B)=\Q(B^{\log})\subset\R$ is
the field generated$/\Q$ by the multiplicities of
$B$, respectively of $B^{\log}$.

\end{description}
Since $\Q(B)$ is countable (small)
it is easy to find required $h_i$ as
small perturbation of multiplicities for
a divisor $H$ with ample $K+B^{\log}+H$.

 \step4
If $K+B_0$ is nef, then $\lambda_0=0$, and
we are done: $(X/Z,B_0=B^{\log})$ is
a strictly log minimal model.

Otherwise there exists
$$
0<\lambda_1=\min\{\lambda\mid
K+B_\lambda \text{ is nef}/Z\}.
$$
$(X/Z,B_{\lambda_1})$ is
a strictly log minimal model too.

\step5 {\em $H$-ordered flips.\/}
As in the proof of Proposition~\ref{resmod},
construction terminates on
the level $\lambda_1$ by a Mori log fibration, or
one can find a log flip (possibly a divisorial contraction)
$X_1\dashrightarrow X_2/Z$ of level $\lambda_1$
with respect to $K+B_0$.
For existence of a Mori log fibration or
of a flipping contraction, one can use
Corollary~\ref{stabray}.
The flip exists by \cite[Theorem~1.1]{HMc} or induction
of Corollary~\ref{exfl}.
By Step~2 the log flop $X_1\dashrightarrow X_2/Z$
with respect to $K+B_{\lambda_1}$ does not
touch $\LCS(X,B_{\lambda_1})=\LCS(X,B_0)$
(see property (1) in Step~3).
Thus it preserves the strictly log minimal model
property of $(X/Z,B_{\lambda_1})$, that is,
$(X^+/Z,B_{\lambda_1})$ is
a strictly log minimal model too.
By Corollaries~\ref{stabray},
\ref{stabwlc} and Addenda~\ref{adstab1}-\ref{adstab2},
as in the proof of Proposition~\ref{resmod},
we obtain an $H$-ordered sequence of extremal log flips
$X_i\dashrightarrow X_{i+1}/Z$ which are
disjoint from $\LCS(X_i,B_{\lambda_i})=\LCS(X,B_0)$.

 \step6 {\em Termination of log flops.\/}
For each level $\lambda>0$, the log flips are
log flops with respect to $K+B_\lambda$.
We claim that there exists at most one such flop,
or equivalently, extremal contractible ray $R$ with
$(K+B_\lambda,R)=0$.

Indeed, let $C$ be a curve$/Z$ in $R$.
Hence $(K+B^{\log}+\lambda H,C)=0$.
Let $C'$ be another curve$/Z$ with
$(K+B^{\log}+\lambda H,C')=0$.
We will verify that $C'$ is also in $R$.
By definition of $\Q(B)$ we have two relations:
$$
\lambda(H,C)=a \text{ and } \lambda(H,C')=a',
$$
with real numbers $a,a'\in \Q(B)$.
Moreover, by (3-4) $a$ and $a'\not=0$.
Otherwise, if $a=0$, all $(H_i,C)=0$ and
$C\equiv 0/Z$; the same holds for $a'$ and $C'$.
Therefore
$$
(H,C)=\frac a\lambda,
(H,C')=\frac {a'}\lambda
$$
and
$$
(H,a'C-aC')=\frac {a'a}\lambda-\frac{aa'}\lambda=0.
$$
Thus, if all $(H_i,a'C-aC')=0$,
by (3) $C'\equiv a'C/a/Z$ and $C'$ in $R$;
$a'/a>0$ by projectivity of $X/Z$.
Otherwise, the $1$-cycle $a'C-aC'$ gives
a nontrivial relation$/\Q(B)$ that contradicts (4).

 \step7 {\em Stabilization of models.\/}
The levels stabilize by Corollary~\ref{equivmod}
for $(X_i/Z,B_i)=(X/Z,B_{\lambda_i})$,
that is, there exists only finitely many
levels: $\lambda_1\ge \lambda_2\ge \dots
\lambda_n>0$.
Of course, we use here a birationally changing model $X$
that means that actually $X$ depends on $\lambda$.
After finitely many log flips we can suppose
condition (2) of Corollary~\ref{equivmod}.
The second statement in (2) holds by Step~2.
Other conditions (1) and (3-4) of Corollary~\ref{equivmod} follow
immediate from construction with
$B=B_{\lambda_1}$ and $B'=B_0$ in (3), and
with $X_i=X_1=X$ in (4)
of the minimal model $(X/Z,B_{\lambda_1})$.
Notice also that models with distinct
levels are nonequivalent by Lemma~\ref{conv}.
Indeed, by construction of $\lambda_{i+1}>0$ and
Step~5 (see also Corollary~\ref{stabwlc} and
its Addendum~\ref{stabwlcad}),
the model $(X/Z,B_{\lambda_{i+1}})$ has
an extremal ray $R$ with $(K+B_{\lambda_{i+1}},R)=0$
different from that of on $(X/Z,B_{\lambda_i})$:
$(K+B_{\lambda_i},R)>0$, and
thus the model is not equivalent to $(X/Z,B_{\lambda_j})$
for $j\le i$.

Finally, the last statement in the theorem follows
from two facts: the numerical log Kodaira dimension
of each minimal model $(X/Z,B)$ is $\ge 0$, and
this is equivalent to the pseudo-effective property of
$K+B$ (cf. Corollary~\ref{closed}).

\end{proof}

\begin{proof}[Proof of Revised Reduction]
The main idea of the reduction is to find
an ordered sequence
of pl flips which, for any level,
has at most one nonspecial flip
of this level on each reduced component of the boundary.
Therefore termination of those flips amounts to
special termination and stabilization of models.
Perturbing and subtracting $H$,
we can achieve this.

We use the proof-construction and
notation from \cite[\S 4]{ISh}.
Construction of
strictly log minimal model $(\oV/Y,B_{\oV}^{\log}+H_{\oV})$
uses only pl flips, and
special termination onto lc centers of codimension $2$
for which log termination in dimension $d-1$
is enough.

We use the following properties of $H$,
an effective reduced (having only
multiplicities $1$) Cartier divisor on $Y$:
\begin{description}

 \item{\rm (1)}
for any birational contraction $\tau\colon Y'\to Y$
with normal $\Q$-factorial variety $Y'$ such that
the prime components of $\tau^*H$ are
all exceptional divisors $E_i$ of $Y'/Y$ and
the proper transforms $H_{Y'}$ of prime components of $H$,
the components and the transform allow to present any numerical class
of divisors$/Y$ \cite[b) p.~65]{ISh}; and

 \item{\rm (2)}
there exists a linear numerical relation between
the components of $\tau^*H$:
$$
\tau^*H=H_{Y'}+\sum a_i E_i=
\sum H_i+\sum a_i E_i\equiv 0/Y
$$
with positive integral $a_i$;
and

 \item{\rm (3)}
the support of $\tau^*H$ is
the reduced part of $B_{Y'}^{\log}+H_{Y'}$ \cite[a) p.~65]{ISh}.
\end{description}
Notice that in both appearances in (2-3), that is,
in the relation and in the boundary,
$H_{Y'}$ is reduced.
All models in our constructions
satisfies the assumption in (1) and
thus properties (1-3).
(Moreover, in what follows,
$\tau$ is an isomorphism over $Y\setminus \Supp H$.)

The next part, to subtract $H_{\oV}$, is quite different.
It should be modified
to use only special termination and
terminal termination in dimension $d$.
We can suppose as usually that $B$ and
$B_{\oV}^{\log}$ are $\Q$-boundaries
(or, in all integral independencies below,
independence should be with multiplicities
of $B$ as well; cf. Step~3, (4) in
the proof of Theorem~\ref{mod}).

\step{1} {\em Perturbation of $H$.\/}
There exists a log flop (not elementary)
$(\oV/Y,B_{\oV}^{\log}+H_{\oV})$ and
a boundary $\Gamma H_{\oV}=\sum \gamma_i H_i$
such that:
\begin{description}

 \item{\rm (4)\/}
$(\oV/Y, B_{\oV}^{\log}+\Gamma H_{\oV})$ is
a strictly log minimal model; and

 \item{\rm (5)\/}
the multiplicities $\gamma_i$ are integrally
(or rationally) independent:
$$
\sum n_i \gamma_i=n, \text{ all } n_i, n\in \Z\Longrightarrow
\text{ all } n_i=0, n=0;
$$

\end{description}
in particular, each $0<\gamma_i<1$.
Moreover, we need $\gamma_i$ arbitrary close to $1$:
$0\ll \gamma_i<1$.
Equivalently, there exists an effective divisor
$\Delta=\sum \delta_i H_i$ with $0<\delta_i\ll 1$
such that:
\begin{description}

 \item{\rm (4)'\/}
$(\oV/Y, B_{\oV}^{\log}+H_{\oV}-\Delta)$ is
a strictly log minimal model; and

 \item{\rm (5)'\/}
the multiplicities $\delta_i$ are integrally
independent.

\end{description}
Indeed,then we can take $\gamma_i=1-\delta_i$.
Actually, it is enough to construct such
a model$/W$:
\begin{description}

 \item{\rm (4)''\/}
$(\oV/W, B_{\oV}^{\log}+H_{\oV}-\Delta)$ is
a strictly log minimal model;

\end{description}
where $W/Y$ is a lc model of
$(\oV/Y,B_{\oV}^{\log}+H_{\oV})$, that is,
the contraction $\oV\to W$ is given
by $K_{\oV}+B_{\oV}^{\log}+H_{\oV}$.
This and other similar
models exist by the LSEPD trick
\cite[10.5]{Sh92} and Corollary~\ref{sample}.
By construction $K_{\oV}+B_{\oV}^{\log}+H_{\oV}$
is ample on $W/Y$, $\equiv 0/W$, and
$-\Delta$ is nef and actually semiample
on $\oV/W$.
Hence, for any $0<\delta\ll 1$,
$K_{\oV}+B_{\oV}^{\log}+H_{\oV}-\delta \Delta$ is
nef on $\oV/Y$.
This gives a model in (4)' with $\Delta:=\delta \Delta$.
Indeed, $K_{\oV}+B_{\oV}^{\log}+H_{\oV}-\delta\Delta$
is nef$/Y$  by construction.
On the other hand, by construction
$(\oV,B_{\oV}^{\log}+H_{\oV})$ is lc,
and by (4)'' and monotonicity  \cite[1.3.3]{Sh92}
$(\oV,B_{\oV}^{\log})$ is
strictly log terminal.
Since $\Supp \Delta=\Supp H_{\oV}$,
again by monotonicity  \cite[1.3.3]{Sh92}
$(\oV, B_{\oV}^{\log}+H_{\oV}-\delta\Delta)$
is strictly log terminal.
The multiplicities of $\delta\Delta$ are
$\delta\delta_i$.
They are integrally independent if $\delta\in \Q$.
Therefore (5)' holds.

Construction of a model in (4)''
uses induction on the number of
irreducible components of $H_{\oV}$.
If it is zero we are done (see Step~5 below).
We start by subtracting the first component
$D=H_1$ of $H_{\oV}$.
This gives new strictly log minimal model
$(\oV/W,B_{\oV}^{\log}+H_{\oV}-D)$.
Since $K_{\oV}+B_{\oV}^{\log}+H_{\oV}\equiv 0/W$,
we need to consider only log flips in
the extremal rays $R$ of $\cNE(\oV/W)$ with
$(D,R)>0$, or
equivalently, to construct a log minimal model
for $(\oV,B_{\oV}^{\log}+H_{\oV}-D)$.
As in the usual reduction the flips are
pl in dimension $d+1$, and thus exist by our assumptions.
The only problem is termination.
Note also that in a birational$/Y$ situation
a Mori log fibration is impossible.
By special termination and
log termination in dimension $d-1$,
a problem is termination in dimension $d$.
More precisely, it is enough to consider
in addition the case of log flips on
a reduced component $F\not=D$ of $B_{\oV}^{\log}+H_{\oV}-D$
with $(F,R)<0$.
Such a component exists by properties (2-3).
According to special termination,
we can consider only flips outside
of the reduced components of the adjoint boundary $B_F$
for $B_{\oV}^{\log}+H_{\oV}-D$ on $F$
which consists of intersections with
components of $H_{\oV}-D-F+\sum E_i$
(see the log adjunction \cite[3.2.3]{Sh92}),
where $\sum E_i$ is the reduced part of $B_{\oV}^{\log}$.
Let $C$ be a flipped curve of one of those log flips,
that is, $C\cap\Supp (H_{\oV}-D-F+\sum E_i)=\emptyset$, and
$(D,C)<0$.
Let $C'$ be the next flipping curve.
Then $(D,C')>0$.
By construction the support of $H_{\oV}$
and exceptional divisors$/Y$
allow to present any numerical class
of divisors$/Y$ (see (1) and (3) above), in particular,
so does any numerically ample divisor$/Y$.
Thus there exists an ample linear combination $aF+bD$
on $C$ and $C'$ with real multiplicities
$a,b$.
On the other hand there exists a nontrivial
(all coefficients are positive) linear relation
between the support of $H_{\oV}$
and exceptional divisors $E_i/Y$ (see (2-3) above).
Hence near $C$ and $C'$ the ample divisor
is given by $cD$ with a real number $c$ which
is impossible under inequalities $(D,C)<0$ and $(D,C')>0$.
Therefore two subsequent flips of this
type do not (co)exist and we get the termination.
Thus we constructed a log pair
$(\oV/W,B_{\oV}^{\log}+H_{\oV}-D)$
which is a strictly log minimal model over $W$, and
let $W_1/W$ denotes its lc model.
Note that each model
$(\oV,B_{\oV}^{\log}+H_{\oV}-D)$
is strictly log terminal because
the starting model
$(\oV,B_{\oV}^{\log}+H_{\oV}-D)$ is
strictly log terminal, and
the log flips preserve this.

Now by induction we can suppose that
a strictly log minimal model
$(\oV/W_1,B_{\oV}^{\log}+H_{\oV}-D-\Delta')$
is constructed where
$\Delta'=\sum_{i\not=1} \delta_i' H_i$, and
the multiplicities $0<\delta_i'\ll 1$ are
integrally independent.
Since $K_{\oV}+B_{\oV}^{\log}+H_{\oV}-D\equiv 0/W_1$ and
$K_{\oV}+B_{\oV}^{\log}+H_{\oV}\equiv 0/W$,
by construction $D=H_1\equiv 0/W_1$.
Hence (1-3) holds on all models $Y'/W_1$ without $H_1$
which are obtained from $\oV/W_1$ by log flops
with respect to $K_{\oV}+B_{\oV}^{\log}+H_{\oV}-D\equiv 0/W_1$
because $H_1\equiv 0/W_1$ on them too
by the assumption $K_{\oV}+B_{\oV}^{\log}+H_{\oV}\equiv 0/W$.
This means that all numerical classes and relations
are considered over $W_1$.
In (3) we can replace $\tau^*H$ by $\tau^*H-H_1$, and
the starting model $(\oV,B_{\oV}^{\log}+H_{\oV}-D)$
is strictly log terminal.
Therefore a required model $\oV/W_1$ exists.
By construction $-D$ is ample on $W_1/W$, and
$-D-\Delta'$ is nef and actually semiample
on $\oV/W_1$.
Hence, for any $0<\delta'\ll \delta\ll 1$,
$-\Delta=-\delta D-\delta'(D+\Delta')$ is
nef on $\oV/W$.
This gives a model in (4)''.
Indeed, $K_{\oV}+B_{\oV}^{\log}+H_{\oV}-\Delta$
is nef$/W$ because
$K_{\oV}+B_{\oV}^{\log}+H_{\oV}\equiv 0/W$.
On the other hand, by construction
$(\oV,B_{\oV}^{\log}+H_{\oV})$ is lc,
and by monotonicity  \cite[1.3.3]{Sh92}
$(\oV,B_{\oV}^{\log})$ is
strictly log terminal.
Since $\Supp \Delta=\Supp H_{\oV}$,
again by monotonicity  \cite[1.3.3]{Sh92}
$(\oV, B_{\oV}^{\log}+H_{\oV}-\Delta)$
is strictly log terminal.
The multiplicities of $\Delta$ are
$\delta_1=\delta+\delta'$ and $\delta_i=\delta'\delta_i',
i\not=1$.
We can take
integrally independent $\delta, \delta_i',i\not=1$,
and $\delta'\in\Q$, .
Then (5)' holds: $\delta_i$ are integrally independent.

 \step2 {\em $H$-ordered flips.\/}
Set $B'=B_{\oV}^{\log}+\Gamma H_{\oV}-
\lambda_{\max} H_{\oV}$
where $\lambda_{\max}$ is the maximal number $\lambda$
such that $B_{\oV}^{\log}+\Gamma H_{\oV}-\lambda H_{\oV}$
is a boundary.
It easy to find that
$\lambda_{\max}=\min\{\gamma_i\}$, and
if say $\lambda_{\max}=\gamma_1$, then
$B'=B_{\oV}^{\log}+\sum_{i\not=1}(\gamma_i-\gamma_1)H_i$.
By construction and monotonicity  \cite[1.3.3]{Sh92}
$(\oV,B')$ and $(\oV,B_{\oV}^{\log}+
\Gamma H_{\oV}=B'+\lambda_{\max}H_{\oV})$
are strictly log terminal with reduced components
$E_i$, exceptional$/Y$, and the second model is log minimal.
Note also that by (5):
\begin{description}

 \item{\rm (5)''\/}
the multiplicities $\gamma_i-\gamma_1,i\not=1$, are integrally
independent, and $0<\gamma_i-\gamma_1\ll 1$.

\end{description}
To construct a log flip in the theorem,
we will find a strictly log minimal model
of $(\oV/Y,B')$ subtracting $H_{\oV}$ from
$B'+\lambda_{\max}H_{\oV}$ as in Proposition~\ref{resmod}.
Existence of $H_{\oV}$-ordered log flips
in this situation is more straightforward.
Indeed, for each level $\lambda\le \lambda_{\max}$,
we can use log flops over a lc model $W/Y$ of
$(\oV/Y,B'+\lambda H_{\oV})$, and
the fact that $\oV/Y$ is FT$/W$.
After construction of a log minimal model$/W$ we
can convert it into a log minimal model$/Y$ as in Step~1, or
use Corollary~\ref{stabwlc} with LSEPD trick.
Note also that a Mori log fibration is impossible.

On the other hand, each flip is pl as in usual reduction.
If it is a log flip in an extremal ray $R$
then by construction
$(K_{\oV}+B'+\lambda H_{\oV},R)=0$
and $(H_{\oV},R)>0$.
Hence by (2) there exists $E_i$ with $(E_i,R)<0$, and
$E_i$ is a reduced component of
$K_{\oV}+B'$.
Since $\lambda>0$,
$(\oV,B')$ is strictly log terminal.

 \step3 {\em Termination of log flops.\/}
For each level $\lambda>0$, the log flips are
log flops with respect to $K_{\oV}+B'+\lambda H_{\oV}$.
By special termination of flips, we can
suppose that after finitely many steps
all the next flips are {\em nonspecial\/} on $\cup E_j$, that is,
log flips in extremal rays $R$ such that, for
any curve $C/Y$ of $R$,
$C$ intersects only one reduced component, say $E_1$.
Actually, $C\subset E_1$, and
$(K_{\oV}+B'+\lambda H_{\oV},C)=0$.
We claim that there exists only one such extremal ray
for $E_1$.
Therefore such log flips terminate.

Let $C'\subset E_1$ be a curve$/Y$ with
$(K_{\oV}+B'+\lambda H_{\oV},C')=0$ and
disjoint from $E_i$ with $i\not=1$.
Since $K_{\oV}+B_{\oV}^{\log}$ is a $\Q$-divisor,
then both relations can be
transform into relations for $\lambda$ and
multiplicities of $B'$:
$$
\lambda (H_{\oV},C)+\sum_{i\not=1}(\gamma_i-\gamma_1)
(H_i,C)=r
$$
and
$$
\lambda (H_{\oV},C')+\sum_{i\not=1}(\gamma_i-\gamma_1)
(H_i,C')=r'
$$
where $r,r'\in\Q$.
Note that $(H_{\oV},C)$ and $(H_{\oV},C')\not=0$
because otherwise we get a rational relation
which contradicts (5)''.
Indeed, if $(H_{\oV},C)=0$ and
all intersections $(H_i,C)=0,i\not=1,$
then $(H_1,C)=0$ too since $H_{\oV}=\sum H_i$.
This is impossible by (1) and (3) because by construction
$(E_i,C)=0, i\not=1$, and by (2) $(E_1,C)=0$.
The same holds for $C'$.
Similarly we verify that $C\equiv c C'/Y$
with $c=(H_{\oV},C)/(H_{\oV},C')\not=0$, e.g.,
$(H_{\oV},C)=c(H_{\oV},C')$.

To do this we eliminate $\lambda$ and
obtain one relation:
$$
(H_{\oV},C')(\sum_{i\not=1}(\gamma_i-\gamma_1)(H_i,C))-
(H_{\oV},C)(\sum_{i\not=1}(\gamma_i-\gamma_1)(H_i,C'))
=r''
$$
with $r''\in \Q$.
By (5)'' this is possible only if
$$
(H_{\oV},C')(\gamma_i-\gamma_1)(H_i,C)=
(H_{\oV},C)(\gamma_i-\gamma_1)(H_i,C').
$$
Again by (5)'' each $\gamma_i-\gamma_1\not=0$
for $i\not=1$.
Hence $(H_i,C)=c(H_i,C')$ for $i\not=1$.
On the other hand, this implies that
$(H_1,C)=c(H_1,C')$ because $H_{\oV}=\sum H_i$.
Then we can prove that $(E_i,C)=c(E_i,C')$
for all $E_i$.
Therefore by (1) $C\equiv cC'/Y$, and
we have the only possibility on $E_1$
for a nonspecial log flip of level $\lambda$.

 \step4 {\em Stabilization of models.\/}
This means that the levels $\lambda_{\max}\ge \lambda_1
\ge \dots >0$ are stabilizes.
This follows from Corollary~\ref{equivmod}.
By special and divisorial termination, after finitely many
log flips, we can suppose only nonspecial log flips
which actually are log flips on corresponding
reduced component $F=E_j$.
Thus conditions (1--2) of Corollary~\ref{equivmod}
hold for $(X_i/Z,B_i)=(F/Y,B_F)$ where
$B_F$ is the adjoint boundary for $B'+\lambda_i H_{\oV}$
on $F=E_j$.
For fixed $F=E_j$, we consider only corresponding
levels $\lambda_i$ (truncation) and models
$(X_i/Z,B_i)$.
Then condition (3) of Corollary~\ref{equivmod}
holds for the adjoint boundary $B$ on $X_i=F$
of pair $(\oV,B')$, and respectively so does
the adjoint boundary $B'$ of $(\oV,B'+\lambda_1 H_{\oV})$.
By Step~2 each $(\oV/Y,B'+\lambda_i H_{\oV})$
is a strictly log minimal model.
Hence construction and adjunction
give (4) of Corollary~\ref{equivmod}
with $X_i$ corresponding to $\lambda_1$.
The models of $F/Y$ are not equivalent
that implies stabilization: there are
only finitely many levels $\lambda_i$.
Indeed, models with distinct
levels are nonequivalent by Lemma~\ref{conv}.
By definition of $\lambda_{i+1}>0$ and
Step~2 (see also Corollary~\ref{stabwlc} and
its Addendum~\ref{stabwlcad}),
the model $(\oV/Y,B'+\lambda_{i+1} H_{\oV})$ has
an extremal ray $R$ with
$(K_{\oV}+B'+\lambda_{i+1} H_{\oV},R)=0$
different from that of on $(\oV/Y,B'+\lambda_i H_{\oV})$:
$(K_{\oV}+B'+\lambda_i H_{\oV},R)>0$,
and the ray is supported on $F=E_j$.
Thus the model $(X_{i+1}/Z,B_{i+1})$
is not equivalent to $(X_i/Z,B_i)$.

 \step5 {\em Flip.\/}
We claim that the cn model $(X^+/Y,B^+)$
of $(\oV/Y,B'+\lambda H_{\oV})$
with $0<\lambda\ll 1$ is a required log flip.
A contraction to the cn model exists
again by LSEPD trick and Corollary~\ref{sample}.
Since the multiplicities $\gamma_i-\gamma_1$ of
$B'+\lambda H_{\oV}=
B_{\oV}^{\log}+\lambda H_1+
\sum_{i\not=1}(\lambda+\gamma_i-\gamma_1) H_i$ are
small and $X$ is $\Q$-factorial,
the boundary $B+\lambda H_1+
\sum_{i\not=1}(\lambda+\gamma_i-\gamma_1) H_i$
on $X$ is klt and $K+B+\lambda H_1+
\sum_{i\not=1}(\lambda+\gamma_i-\gamma_1) H_i$
is numerically negative on $X/Y$.
Notice also that $(\oV/Y,B'+\lambda H_{\oV})$
is a strictly log minimal model
of $(X/Y,B+\lambda H_1+
\sum_{i\not=1}(\lambda+\gamma_i-\gamma_1) H_i)$.
Thus by monotonicity \cite[Lemma~2.4]{ISh},
$B'+\lambda H_{\oV}$ does not have reduced components and by (3)
the contraction $\oV\to Y$ is small.
Thus the cn model is the log flip of
$(X/Y,B+\lambda H_1+
\sum_{i\not=1}(\lambda+\gamma_i-\gamma_1) H_i)$.
Since the contraction $X\to Y$ is extremal,
its flip is unique \cite[Corollary~3.6]{Sh00}, and $(X^+/Y,B^+)$
is a required log flip.

\end{proof}

\begin{proof}[Proof of Revised Induction]
By methods of \cite{Sh00} \cite{HMc}
it is sufficient to establish \cite[Theorem~7.2]{HMc}
in dimension $d$.
It can be obtained from Theorems~\ref{mod} and \ref{fin}:
the former gives existence of models $W_i$ and
the later gives the finiteness of models.
Moreover, we can assume that $(X/Z,\Delta)$
is a strictly log terminal pair with
a birational contraction $X/Z$.
Thus we need both theorem only in
the birational situation, and
needed terminal termination is also
only for birational pairs.
Notice that for klt pairs their $\Q$-factorialization
can be obtain by Theorem~\ref{mod} as a strictly log minimal model
with boundary multiplicities $1$ for
the exceptional divisors;
actually, it is enough existence of log flips in
dimension $d$ and special termination.
The log flips exist by \cite[Theorem~1.1]{HMc} or induction
in Corollary~\ref{exfl}.

Using LSEPD trick we can slightly increase $\Delta$
and assume that $\Delta\in V\subseteq \sD_F$
for $F=\Supp \Delta$; and $F$ contains all
exceptional divisors of $X/Z$.
Any $\Theta\in \sD_F$ sufficiently close to $\Delta$
is an $\R$-boundary.
Thus $(X/Z,\Theta)$ has a strictly log minimal model
$(W_i/Z,\psi_{i*}\Theta)$ where
$\psi_i\colon X\dashrightarrow W_i$ is a birational rational $1$-contraction.
It is $1$-contraction by definition and
monotonicity \cite[Lemma~2.4]{ISh}.
Since $X/Z$ is birational, Mori log fibrations are
impossible.
(According to the proof of Theorem~\ref{mod}
after Step~3, we can decompose $\psi_i$ into
a composition of log flips that gives (1) in
\cite[Theorem~7.2]{HMc};
cf. also Proposition~\ref{resmod}.
However it is not important for pl flips,
in particular, for \cite[Corollary~7.3]{HMc}.)
By construction (2--3) of \cite[Theorem~7.2]{HMc} hold,
and (4) see in the proof of \cite[Theorem~7.2]{HMc}.

Thus we need to establish the finiteness up
to equivalence of models.
If this does not hold we have a convergent
sequence of $\R$-boundaries $\Theta_i\in \sD_F$:
$\lim_{i\to\infty}\Theta_i=\Delta$,
with pairwise nonequivalent wlc models
$(X_i/Z,B_i)=(W_i/Z,\psi_{i*}\Theta_i)$.
This contradicts Theorem~\ref{fin}
(cf. the proof of Corollary~\ref{equivmod}).
For a subsequence of models,
assumption (1) of Theorem~\ref{fin} holds
by construction and the birational property of $X/Z$
(still we need to construct $(X/Z,B)$ and
verify (1) for it).
We construct an appropriate model $(X/Z,B)$ as
a modification of a strictly log minimal model
$(W/Z,\Delta)$.
The model is birationally larger that $X_i$ for
a subsequence with $i\gg 0$:
$W\dashrightarrow X_i$ is a birational rational $1$-contraction.
Indeed, we can construct it from $X$ by
a sequence of ordered log flips of $(X/Z,\Delta)$
(see the proof of Theorem~\ref{mod}).
By stability of negative contractions and
Main Lemma we never contract a component of $\Delta$
which are not contracted on $X_i$ with all $i\gg 0$
(cf. Step 3 in the proof of Theorem~\ref{fin}).
Since the sequence of log flips is finite
we get required subsequence of models $(X_i/Z,B_i)$.
Hence the model $(W/Z,\Delta)$ is between
$(X/Z,B)$ and $(\oX/Z,B_{\oX})$ in the proof of
Theorem~\ref{fin}, and we can construct $(\oX/Z,B_{\oX})$
as in the proof.
Otherwise we consider a lc model $(W/Z,\Delta)$
which exists by Corollary~\ref{sample}.
It contracts all prime divisors $D$ on $X$ which
are exceptional on varieties $X_i$ with $i\gg 0$.
Otherwise $K+\Delta$ is big on $D/Z$, and
so does $K+\Theta_i$ with $i\gg 0$ which
contradicts construction:
each log flip or divisorial contraction
of $(X/Z,\Theta_i)$
preserves this property \cite[Proposition~3.20]{Sh00}.
Thus we can construct $(X/Z,B)$ as
a crepant blowup of $(W/Z,\Delta)$ with
blown up divisors $D$ the same as for models $X_i/Z$.
This model also can be constructed as
a strictly log minimal model with multiplicities
$1$ for other exceptional divisors$/W$ on
a starting model.
Pair $(X/Z,B)$ satisfies (1).

Assumption (2) is void because $(X,B)$ and
$(X,B_i)$ are klt.

In (3) we take all boundary components of $B$ and
the exceptional divisors $E$ with $b(E,X,B)\ge 0$.
Then assumptions (4--5) hold for an appropriate
subsequence as in Step~2 in the proof of Corollary~\ref{equivmod}.
\end{proof}

\bigskip
\noindent Department of Mathematics,\\ Johns Hopkins University,\\
Baltimore, MD--21218, USA\\ e-mail: shokurov@math.jhu.edu

\noindent
Steklov Mathematical Institute,\\ Russian Academy of Sciences,\\
Gubkina str. 8, 119991, Moscow, Russia

\end{document}